\documentclass[journal,twoside,print]{ieeecolor_arxiv}
\usepackage{generic}
\usepackage{cite}
\usepackage{amsmath,amssymb,amsfonts}
\usepackage{algorithmic}
\usepackage{graphicx}
\usepackage{textcomp}
\newcommand{\TitleNamemain}{Sensitivity-Conditioning: \\ Beyond Singular Perturbation for Control Design on Multiple Time Scales}

\def\BibTeX{{\rm B\kern-.05em{\sc i\kern-.025em b}\kern-.08em
    T\kern-.1667em\lower.7ex\hbox{E}\kern-.125emX}}

\usepackage{graphicx}      

\usepackage{commath}
\usepackage{mathtools}
\usepackage{steinmetz}		
\usepackage{supertabular}	
\usepackage{tabularx}		
\usepackage{algorithm}		
\usepackage{algorithmic}		
\usepackage{amsmath,amssymb}

\usepackage{amsthm} 
\usepackage{multirow}
\usepackage{caption}
\usepackage{subcaption}
\usepackage{adjustbox}

\newtheorem{rem}{Remark}
\newtheorem{defi}{Definition}
\newtheorem{prop}{Proposition}
\newtheorem{thm}{Theorem}
\newtheorem{ass}{Assumption}
\newtheorem{lemma}{Lemma}
\newtheorem{cor}{Corollary}
\usepackage[usenames,dvipsnames]{xcolor}
\usepackage{soul}
\usepackage{url} 
\usepackage{afterpage}
\let\labelindent\relax
\usepackage{enumitem}
\usepackage{cancel}

\usepackage{tikz}
\usetikzlibrary{shapes,arrows}
\usepackage[american]{circuitikz}

\usepackage{tikz}
\usetikzlibrary{calc}
\newcommand{\hcancel}[5]{%
    \tikz[baseline=(tocancel.base)]{
        \node[inner sep=0pt,outer sep=0pt] (tocancel) {#1};
        \draw[black] ($(tocancel.south west)+(#2,#3)$) -- ($(tocancel.north east)+(#4,#5)$);
    }%
}%

\tikzset{
block/.style = {draw, fill=white, rectangle, minimum height=1em, minimum width=1em},
block2/.style = {draw, fill=white, rectangle, minimum height=3em, minimum width=3em},
tmp/.style  = {coordinate}, 
sum/.style= {draw, fill=white, circle, node distance=1cm},
input/.style = {coordinate},
output/.style= {coordinate},
pinstyle/.style = {pin edge={to-,thin,black}
}
}

\newcommand{\normsz}[1]{\lVert #1 \rVert}

\newcommand{\rr}[1]{{\color{black}#1}}

\newcommand{\w}[1]{\textcolor{white}{#1}}

\newcommand{\saverio}[1]{{\color{black}#1}}

\begin{document}

\title{\TitleNamemain}

\author{Miguel Picallo, Saverio Bolognani, Florian D{\"o}rfler 
\thanks{Funding by the Swiss Federal Office of Energy through the project “Renewable Management and Real-Time Control Platform (ReMaP)” (SI/501810-01) and the ETH Foundation is gratefully acknowledged.}
\thanks{The authors are with the Automatic Control Laboratory at ETH Z{\"u}rich, Switzerland.
(emails: \{miguelp,bsaverio,dorfler\}@ethz.ch)}%
}

\maketitle

\begin{abstract}
A classical approach to design controllers for interconnected systems is to assume that the different subsystems operate at different time scales, then design simpler controllers within each time scale, and finally certify stability of the interconnected system via singular perturbation analysis. In this work, we propose an alternative approach that also allows to design the controllers of the individual subsystems separately. However, instead of requiring a sufficiently large time-scale separation, our approach consists of adding a feed-forward term to modify the dynamics of faster systems in order to anticipate the dynamics of slower ones. We present several examples in bilevel optimization and cascade control design, where our approach improves the performance of currently available methods.
\end{abstract}


\IEEEpeerreviewmaketitle

\section{Introduction}
Interconnected and nested systems are ubiquitous in control applications
, but they may be challenging to analyse and design. If interconnected systems are composed by subsystems operating on multiple time scales \cite{chow1982timescalesint} \saverio{and a normal hyperbolicity condition holds \cite{kuehn2015multiple}}, then each time scale can be studied independently, substituting dynamics of faster time scales by algebraic equations \cite{tikhonov1952systems}. Such systems appear in engineering applications like power systems \cite{peponides1982singular, subotic2019lyapunov}, biological systems \cite{muratori1992foodchain}, motion control \cite{tavasoli2007motion2timescale}, electrical drives \cite{mezouar2007elecdrive2timesca}, etc. In that context, time-scale separation arguments, like singular perturbation analysis \cite{kokotovic1976singpertoverview, kokotovic1999singular}, allow to certify when the stability guarantees derived in each separate time scale are preserved in the interconnected system. 
Standard singular perturbation considers only two time scales \cite{khalil2002nonlinear}, although it can be extended to multiple ones \cite{grammel2004nonlinear, kuehn2015multiple,subotic2019lyapunov}.

Besides analysis, singular perturbation is also a powerful tool for control design \cite{kokotovic1984applications}, for example as a model reduction technique \cite{koko1968singpertmodelred}: complex systems on a single time scale can be artificially separated into subsystems on different time scales, and thus simplify their analysis and controller design. Singular perturbation analysis can then provide additional conditions, for example on the control parameters \cite{subotic2019lyapunov}, to ensure that the interconnected system remains stable. Some examples of these applications are hierarchical control architectures, like cascade control \cite{lee1998pidcascade}, or iterative optimization algorithms, like dual ascent \cite{bertsekas1997nonlinear}, interior point methods \cite{boyd2004convex}, etc. However, for more than two time scales such singular perturbation conditions may be hard to derive, unless the interconnection present a specific structure  \cite{grammel2004nonlinear,subotic2019lyapunov}. More importantly, since artificial time-scale separation slows down some subsystems with respect to others, it poses a fundamental limit on the convergence rate of the interconnected system.


\saverio{In this work, we consider interconnected control systems in which the individual subsystems are designed and stabilized (e.g., by means of control) on separate time scales, and we are interested in preserving the overall system stability of the interconnection in a single time scale. Unlike the singular perturbation approach, we propose a single-time scale interconnection that guarantees closed-loop stability} without imposing additional conditions on control parameters, nor slowing down any subsystem with respect to others. Additionally, our approach can deal with general interconnection structures, where the dynamics of each subsystem may depend on the states of all other subsystems. \rr{Our proposed interconnection can be interpreted as a transient feed-forward term in faster systems, that anticipates the dynamics of slower ones. For that, it uses the sensitivity of the fast system's steady state with respect to the slower system's state. Therefore, we term this approach the sensitivity-conditioning. 
}

This new interconnection is inspired by recently proposed optimization algorithms to solve problems that are usually represented on multiple time scales: the prediction-correction algorithms for time-varying optimization \cite{simonetto2016predcorr, fazlyab2018predcorintpoint}, the advance-steps in nonlinear model predictive control \cite{zavala2009advstepnmpc, diehl2009effnmpc}, and the opponent-learning awareness games \cite{foerster2017learning, wang2019followridge}. These algorithms use the nonlinear optimization sensitivity \cite{fiacco1976sensitivity, jittorntrum1984solution} to generate feed-forward terms that improve their convergence. 
\rr{Our approach also relates to classic backstepping \cite[Ch.~14]{khalil2002nonlinear} in the context of overcoming time-scale separation limitations. However, unlike backstepping, our approach does not require to know a stabilizing state feedback law in closed form. Hence, our approach is implementable in cases where such a feedback law is not available.}

\rr{Our contributions are the following: 
First, we divide the problem of designing the interconnection of two subsystems into a design problem of separate time-scales and a conditioning of their interconnection. For the latter, we define the sensitivity-conditioning approach, and we show how it corresponds to an additional transient signal to be exchanged between the two subsystems. Second, we prove that the sensitivity-conditioning approach ensures that the interconnected single-time-scale system has the same local (and even global, under further conditions) exponential stability properties as the multiple time scales system where each subsystem evolves on a different time scale. Third, we show how some degrees of freedom in the proposed design method can be used to improve the convergence rate of the interconnected system, and we provide robustness guarantees with respect to model errors. Fourth, we demonstrate the applicability of our approach on two control design problems: cascade control \cite{lee1998pidcascade} and bilevel optimization \cite{bard1998practical}. 
Finally, we show how to extend our approach to multiple time scales.}

The rest of this paper is structured as follows: Section~\ref{sec:moti} presents the type of systems that we consider and motivates our sensitivity-conditioning approach. Section~\ref{sec:predsens2sys} introduces the sensitivity-conditioning for two interconnected systems. Sections~\ref{sec:casccont} and \ref{sec:bilopt} show the applications examples. Section~\ref{sec:multits} shows how the sensitivity-conditioning can be extended to multiple-time-scales systems, and discrete-time systems. Finally, Section~\ref{sec:conc} presents some conclusions.

\section{Motivation}\label{sec:moti}

\saverio{
Consider the interconnection of two systems described by the vector fields $f_i(\cdot)$ on $x_i \in \mathbb{R}^{n_i}$, respectively for $i=1,2$:
\begin{equation}\label{eq:noprecondsys}
\begin{split}
\dot{x}_1 &= f_1(x_1,x_2) \\ 
\dot{x}_2 &= f_2(x_1,x_2).
\end{split}
\end{equation}

The study (or design) of such interconnection is challenging in general. 
One way to tackle these analysis or design problems is to assume that the two subsystems in \eqref{eq:noprecondsys} evolve on separate time scales: $x_2$ evolves on a faster time $\tau$, where $x_1$ is constant, while $x_1$ evolves on a slower time $t$, where the dynamics of $x_2$ are replaced by the algebraic equation $f_2 (x_1,x_2)=0$.
This two-time-scale interconnection is represented with a differential-algebraic-equation system $\Sigma_1$, and a boundary-layer system $\Sigma_2$:
\begin{subequations}\label{eq:sys2ideal0}
    \begin{align}
    \Sigma_1 & : \dot{x}_1 = f_1(x_1,x_2) & \text{s.t. }& f_2(x_1,x_2)=0  \label{seq:sys2ideal0dae} \\
    \Sigma_2 & : \frac{dx_2}{d\tau} = f_2(x_1,x_2) & \text{s.t. }& \frac{dx_1}{d\tau}=0. \label{seq:sys2ideal0bl}
\end{align}
\end{subequations}
Many interconnected systems become simpler to design and control when the two subsystems are assumed to evolve on two separate time scales as in \eqref{eq:sys2ideal0}.
Classical examples are adaptive control \cite{aastrom2013adaptive}, cascade control systems \cite{lee1998pidcascade}, where the state of one system is used as input to the other system, i.e., $x_2$ to control $x_1$, or nested iterative numerical algorithms (e.g., in optimization). 
In the rest of the paper, we assume that the analysis and design of the two-time-scale system \eqref{eq:sys2ideal0} are tractable problems, and we provide some examples of how this is done for specific applications in Sections~\ref{sec:casccont} and \ref{sec:bilopt}.

Clearly, any statements on the steady-state behavior and the stability of the two-time-scale system~\eqref{eq:sys2ideal0} does not automatically hold true for the original single-time-scale system~\eqref{eq:noprecondsys}.
One standard way to ensure that the properties of \eqref{eq:sys2ideal0} extend to \eqref{eq:noprecondsys} is to enforce a sufficient time-scale separation between the two subsystems and then employ the tools of singular perturbation analysis \cite[Ch.~11]{khalil2002nonlinear}.
Under the assumption that $f_2(x_1,\cdot)$ has a finite number of isolated roots $x_2^s(x_1)$, one can define the \textit{standard} singular perturbation {conditioned} system
\begin{equation}\label{eq:sys2singpert}
\begin{split}
\dot{x}_1 & = f_1(x_1,x_2) \\ 
\epsilon \dot{x}_2 & = f_2(x_1,x_2).
\end{split}
\end{equation}
where $0<\epsilon \ll 1$ is a design parameter to enforce the desired level of time-scale separation.
In the \textit{singular limit} $\epsilon \to 0$, \eqref{eq:sys2singpert} becomes a degenerate system by Tikhonov's Theorem \cite{tikhonov1952systems} and reduces to \eqref{eq:sys2ideal0}.
Singular perturbation analysis allows to guarantee that if both systems in \eqref{eq:sys2ideal0} are asymptotically stable, then the conditioned interconnection \eqref{eq:sys2singpert} is also stable (and has the same equilibria) when $\epsilon$ is below a certain threshold $\bar{\epsilon}$ \cite[Thm.~11.3,4]{khalil2002nonlinear}.
An example of a control design targeting a time-scale separated closed-loop system as in \eqref{eq:sys2singpert} is cascaded control, e.g.,
in power electronics control systems, where time-scale separation does not exist naturally, but has to be imposed artificially \cite{subotic2019lyapunov}. This type of conditioning comes at a cost: as the second subsystem cannot be made arbitrarily fast in practice, the design choice of $\epsilon \ll 1$ necessarily slows down the first subsystem and thus limits the convergence rate and deteriorates the performance of the entire interconnection \eqref{eq:sys2singpert}.

In this paper, we propose an alternative conditioning of the interconnected system \eqref{eq:noprecondsys} without this drawback. For that, we define the conditioned interconnected system $\Sigma$ as
\begin{equation}\label{eq:intersys}
\Sigma: \quad 
M(x_1,x_2) \begin{bmatrix} \dot{x}_1\\ \dot{x}_2 \end{bmatrix}
= 
\begin{bmatrix} f_1(x_1,x_2) \\ f_2(x_1,x_2)\end{bmatrix},
\end{equation}
where $M(\cdot) = \left[ \begin{smallmatrix} M_{11} & M_{12} \\ M_{21} & M_{22} \end{smallmatrix} \right] \hspace{-0.1cm} (\cdot)$ is a general non-singular conditioning matrix, i.e., a generalized time constant, which is a design variable to be chosen. Notice that the singular perturbation conditioned \eqref{eq:sys2singpert} is a special case of \eqref{eq:intersys}, with a specific matrix $M=\left[\begin{smallmatrix} I & 0 \\ 0 & \epsilon I \end{smallmatrix}\right]$, where $I$ is the identity matrix. 
Nonetheless, a general $M$ can represent a much larger class of interconnections, see Table~\ref{tab:comparM} for an illustration. For example, it can represent the fully actuated interconnected control system
\begin{equation}\label{eq:intersysu}
\begin{array}{rl}
\dot{x}_1 & = f_1(x_1,x_2) + u_1 \\ 
\dot{x}_2 & = f_2(x_1,x_2) + u_2,
\end{array}
\end{equation}
where the external control inputs $u_i$ are active only the transient dynamics, i.e., $u_i=0$ when $f_i(\cdot)=0$, thus preserving the steady-states of the original system $\eqref{eq:noprecondsys}$ and of the two-time-scale system \eqref{eq:sys2ideal0}. One example of this general conditioning technique \eqref{eq:intersys} is backstepping in cascade systems \cite[Ch.~14]{khalil2002nonlinear}, which we will review in Section~\ref{sec:casccont}. 
Other examples appear in the design of nested gradient algorithms for continuous-time optimization, which we will discuss in Section \ref{sec:bilopt}.

Note that the singular-perturbation conditioned system \eqref{eq:sys2singpert} uses a  \textit{transient} control action $u_2=(\frac{1}{\epsilon}-1)f_2(x_1,x_2)$ (and $u_1=0$) to induce time-scale separation. In this article, we propose an alternative conditioning matrix $M(\cdot)$ in \eqref{eq:intersys} that, in the form \eqref{eq:intersysu}, corresponds to a  derivative-type control action $u_2 = M_{21}f_1(x_1,x_2)$. 
Loosely, we propose that the $x_2$-dynamics are additionally driven by $u_2 \approx \frac{d}{dt}x_2^s(x_1(t))$, where $x_2^s(x_1)$ is the steady state of the boundary layer system \eqref{seq:sys2ideal0bl} parametrized by $x_1$, i.e., $f_2(x_1,x_2^s(x_1))=0$. As a result, under the dynamics \eqref{eq:intersysu} we have $\frac{d}{dt} \normsz{x_2(t) - x_2^s(x_1(t))}_2^2 = 2 (x_2(t) - x_2^s(x_1(t)))^T\big(f_2(x_1(t),x_2(t))+\cancel{u_2 - \frac{d}{dt}x_2^s(x_1(t))}\big)$. Likewise, for \eqref{seq:sys2ideal0bl} we have $\frac{d}{d\tau} \normsz{x_2(\tau) - x_2^s(x_1)}_2^2 = 2 (x_2(\tau) - x_2^s(x_1))^T f_2(x_1,x_2(\tau))$. In other words, if stability of the instantaneous steady-state $x_2^s(x_1)$ of \eqref{seq:sys2ideal0bl} can be inferred by means of a quadratic Lyapunov function, so can be the stability of the trajectory $x_2^s(x_1(t))$ of \eqref{eq:intersysu}. In either case, the stability analysis of the coupled $(x_1,x_2)$ dynamics reduces to that of a cascade system, and no time-scale separation is required. 

However, generally $x_2^s(x_1(t))$ and its derivative $\frac{d}{dt}x_2^s(x_1(t))$ are not available in closed form. In what follows, we show how to construct an implementable surrogate for $u_2 \approx \frac{d}{dt}x_2^s(x_1(t))$, analyze the system stability without requiring $x_2^s(x_1(t))$, and extend the argument to an arbitrary number of subsystems.
}

\newsavebox{\smlmat}
\savebox{\smlmat}{$M(x_1,x_2) = \left[\begin{smallmatrix} M_1 & 0 \\ 0 & M_2 \end{smallmatrix}\right]$}

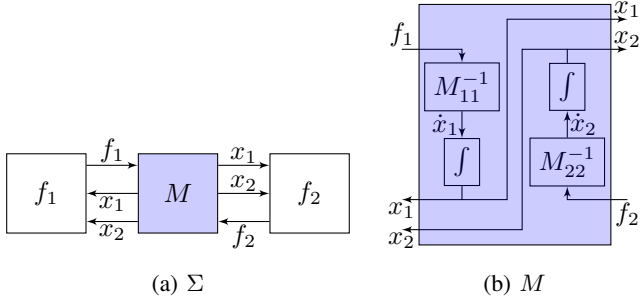
\begin{figure}[t]
     \centering
     \begin{subfigure}[t]{0.28\textwidth}
         \centering
\begin{tikzpicture}[auto, node distance=2cm,>=latex']
    \node [block2, name=system1] (system1) {$f_1$};
    \node [block2, name=system2, right of=system1, node distance = 3.5cm] (system2) {$f_2$};
    \node [block2, name=matrix, right of=system1, node distance = 1.75cm, fill=blue ,fill opacity=.2, draw opacity=1, text opacity = 1] (matrix) {$M$};
    \draw [->] (system1.35) --node[above=-2.5pt,name=f1]{$f_1$} (matrix.145) ;
    \draw [->] (matrix.180) --node[below=-2.5pt,name=dx1]{$x_1$} (system1.0) ;
    \draw [->] (matrix.215) --node[below=-2.5pt,name=dx2]{$x_2$} (system1.325) ;
    \draw [->] (matrix.35) --node[above=-2.5pt,name=dx1r]{$x_1$} (system2.145);
    \draw [->] (matrix.0) --node[above=-2.5pt,name=dx2r]{$x_2$} (system2.180);
    \draw [->] (system2.215) --node[below=-2.5pt,name=f2]{$f_2$} (matrix.325) ;
    
\end{tikzpicture}
    \caption{$\Sigma$}\label{fig:intersys}
     \end{subfigure} %
     \begin{subfigure}[t]{0.2\textwidth}
         \centering
\begin{tikzpicture}[auto, node distance=2cm,>=latex']
    \node [input,name=f1] (f1) {};
    \node [tmp,name=tf1,right of = f1, node distance =0.8cm] (tf1) {};
    \node [block, name=m1, below of = tf1, node distance = 0.5cm] (m1) {$M_{11}^{-1}$};
    \node [block, name=int1, below of = m1, node distance = 1.cm] (int1) {$\int$};
    \node [tmp,name=tx1,below of = int1, node distance =0.5cm] (tx1) {};
    \node [input,name=x1, left of =tx1, node distance = 0.8cm] (x1) {};
    
    \draw[->] (f1) |- node[above=-2.5pt]{$f_1$} (tf1) -| (m1);
    \draw[->] (m1) -- node[left=-2.5pt,pos=0.5,name=dx1]{$\dot{x}_1$}  (int1);
    \draw[->] (int1) |- (tx1) -- node[below=-2.5pt,pos=1]{$x_1$} (x1);

    \node [input,name=f2,right of = x1, node distance = 3cm] (f2) {};
    \node [tmp,name=tf2,left of = f2, node distance =0.8cm] (tf2) {};
    \node [block, name=m2, above of = tf2, node distance = 0.5cm] (m2) {$M_{22}^{-1}$};
    \node [block, name=int2, above of = m2, node distance = 1.cm] (int2) {$\int$};
    \node [tmp,name=tx2,above of = int2, node distance =0.5cm] (tx2) {};
    \node [input,name=x2, right of =tx2, node distance = 0.8cm] (x2) {};
    
    \draw[->] (f2) |- node[below=-2.5pt
    ]{$f_2$} (tf2) -| (m2);
    \draw[->] (m2) -- node[right=-2.5pt,pos=0.5,name=dx2]{$\dot{x}_2$}  (int2);
    \draw[->] (int2) |- (tx2) -- node[above=-2.5pt,pos=1]{$x_2$} (x2);
    
    \node [input,name=x1p,above of = x2, node distance = 0.4cm] (x1p) {};
    \node [tmp, name= tx1p3, right of =tx1, node distance = 0.6cm] (tx1p3) {};
    \node [tmp, name= tx1p4, above of =tx1p3, node distance = 2.4cm] (tx1p4) {};
    \draw[->] (tx1) -- (tx1p3) --(tx1p4) --node[above=-2.5pt,pos=1]{$x_1$} (x1p);
 
    \node [input,name=x2p,below of = x1, node distance = 0.4cm] (x2p) {};
    \node [tmp, name= tx2p3, left of =tx2, node distance = 0.6cm] (tx2p3) {};
    \node [tmp, name= tx2p4, below of =tx2p3, node distance = 2.4cm] (tx2p4) {};
    \draw[->] (tx2) -- (tx2p3) --(tx2p4) --node[below=-2.5pt,pos=1]{$x_2$} (x2p);   
    \node [block,name=bluesq, right of =dx1, node distance = 0.92cm , fill=blue ,fill opacity=.2, draw opacity=1, text opacity = 1, minimum height = 3.2cm, minimum width = 2.55cm] (bluesq) {}; 
    
\end{tikzpicture}
\caption{$M$}\label{fig:diagM}
     \end{subfigure}
        \caption{Block diagram of the conditioned interconnected system $\Sigma$. Fig.~\ref{fig:diagM} shows an example of a constant block-diagonal conditioning matrix~\usebox{\smlmat}.}
        \label{fig:intersysdiagM}
\end{figure}

\section{Sensitivity-conditioning for two systems}\label{sec:predsens2sys}

\subsection{Steady states and sensitivity}

For this section we make the following standard simplifying assumption \cite[Ch.~11]{khalil2002nonlinear}, which we will partially relax later in Section~\ref{sec:multits}:
\begin{ass}\label{ass:sys2sensinv}
The vector fields $f_i(\cdot)$ are continuously differentiable. \rr{For every $x_1$, $f_2(x_1,\cdot)=0$ has a single root $x_2^s$, where the partial derivative $\nabla_{x_2}f_2(x_1,x_2^s)$ is invertible.}
\end{ass}

Under Assumption \ref{ass:sys2sensinv}, the implicit function theorem \cite{krantz2012implicit} guarantees the local existence of a continuously differentiable steady-state map $x_2^s(x_1)$, and gives the sensitivity of this steady state $x_2^s(x_1)$ with respect to $x_1$ as
\begin{equation*}\label{eq:sys2sensss}\begin{array}{l}
     \nabla_{x_1} x_2^s(x_1) \hspace{-0.05cm} = \hspace{-0.05cm} - \hspace{-0.05cm} \nabla_{x_2}f_2(x_1,x_2^s(x_1)\hspace{-0.05cm})^{-1}  \nabla_{x_1}f_2(x_1,x_2^s(x_1)\hspace{-0.05cm}), 
\end{array}
\end{equation*}
where $\nabla_{x_1} x_2^s(x_1) \in \mathbb{R}^{n_2 \times n_1}$. Even though this sensitivity is defined only at points where $x_2=x_2^s(x_1)$, given Assumption~\ref{ass:sys2sensinv} its analytic expression 
is well-defined at any point $x_2$ in a neighborhood of $x_2^s(x_1)$. This allows to define an extended sensitivity 
\begin{equation}\label{eq:sys2sens}\begin{array}{l}
     S_{x_1}^{x_2} (x_1,x_2) := -\nabla_{x_2}f_2(x_1,x_2)^{-1} \nabla_{x_1}f_2(x_1,x_2), 
\end{array}
\end{equation}
which satisfies the restriction $S_{x_1}^{x_2} (x_1,x_2^s(x_1))=\nabla_{x_1} x_2^s(x_1)$.

The steady state map $x_2^s(x_1)$ allows us to redefine the differential-algebraic-equation system \eqref{seq:sys2ideal0dae} as a reduced-order system with reduced vector field $f_1^r(\cdot)$, so that the two-time-scale system \eqref{eq:sys2ideal0} becomes   
\begin{subequations}\label{eq:sys2ideal}\begin{align}
    \Sigma_1&: \dot{x}_1 = f^r_1(x_1) := f_1(x_1,x_2^s(x_1))  \label{seq:sys2idealro}\\
    \Sigma_2&: \frac{dx_2}{d\tau} = f_2(x_1,x_2) \text{ s.t. } \frac{dx_1}{d\tau}=0.
\end{align}
\end{subequations}

Then, each steady state $x_1^s$ satisfying $f^r_1(x_1^s)=0$, defines a steady state $(x_1^s,x_2^s(x_1^s))$ for \eqref{eq:sys2ideal}, and by \cite[Cor.~4.3]{khalil2002nonlinear} it is a locally exponentially stable steady state of the two-time-scale system \eqref{eq:sys2ideal} if and only if $\nabla_{x_2} f_2(x_1^s,x_2^s(x_1^s))$ and
\begin{equation}\label{eq:sys2totder}\begin{array}{l}
     \nabla_{x_1} f^r_1(x_1^s) = \\ \nabla_{x_1}f_1(x_1^s,x_2^s(x_1^s)) + \nabla_{x_2}f_1(x_1^s,x_2^s(x_1)) \nabla_{x_1}x_2^s(x_1^s)
\end{array}
\end{equation}
have eigenvalues with negative real part, and unstable if any of these matrices has any eigenvalue with positive real part.

\begin{rem}\label{rem:stablost}
Both the singular perturbed \eqref{eq:sys2singpert} and the two time-scale systems \eqref{eq:sys2ideal} have the same steady state $(x_1^s,x_2^s(x_1^s))$, but their local exponential stability properties may differ given the value of $\epsilon$, because their Jacobians, and thus their local linearisations, may have different eigenvalues, see Table~\ref{tab:comparM}. 

For example, consider the linear system $f_1(x_1,x_2)=x_1-2x_2$, $f_2(x_1,x_2)=\frac{1}{2}x_1-\frac{1}{2}x_2$, thus $x_2^s(x_1)=x_1$, $\nabla_{x_1} x_2^s(x_1)=1$, $f^r_1(x_1)=-x_1$ and $(x_1^s,x_2^s(x_1^s))=(0,0)$. Then \eqref{eq:sys2ideal} is exponentially stable, since $J_2=\nabla_{x_2}f_2(0,0)=-\frac{1}{2}<0$ and $J_1=\nabla_{x_1} f^r_1(0)=-1<0$. The Jacobian of \eqref{eq:sys2singpert} in Table~\ref{tab:comparM} has indeed negative eigenvalues for $\epsilon<\frac{1}{2}$, but positive ones for $\epsilon>\frac{1}{2}$. Thus \eqref{eq:sys2singpert} is exponentially stable if and only if $\epsilon<\bar{\epsilon}=\frac{1}{2}$. 
\end{rem}


\begin{table*}[t]
\normalsize
  \centering
\caption{Comparison of interconnections $M$}
  \label{tab:comparM}
  \renewcommand{\arraystretch}{1.}
  \begin{tabular}{|c|c|c|c|c|}
     \hline
    Case & Two time scales \eqref{eq:sys2ideal} & Singular perturbed \eqref{eq:sys2singpert} &  Sensitivity-cond. \eqref{eq:sys2predsens} &
    Generalized S-C \eqref{eq:sys2predsensacc}\\ \hline
  
    System $\Sigma$ & {$\begin{array}{rl}
    \dot{x}_1 = & f^r_1(x_1) \\ 
    \frac{dx_2}{d\tau} = & f_2(x_1,x_2)
    \end{array}$} &  
    \multicolumn{3}{c|}{$M \begin{bmatrix} \dot{x}_1 \\ \dot{x}_2 \end{bmatrix} = \begin{bmatrix} f_1(x_1,x_2) \\ f_2(x_1,x_2) \end{bmatrix}$} \\ \hline
  
    \multicolumn{1}{|c|}{\begin{tabular}{c}
         Matrix $M$ 
    \end{tabular}} &
    \multicolumn{1}{m{2.cm}|}{$\begin{bmatrix} I & 0 \\ 0 & 0 \end{bmatrix}$} 
    & \multicolumn{1}{m{2.4cm}|}{$\begin{bmatrix} I & 0 \\ 0 & \epsilon I \end{bmatrix}$} 
    & \multicolumn{1}{m{2.7cm}|}{$\begin{bmatrix} I & 0 \\ -S_{x_1}^{x_2} & I \end{bmatrix}$}
    & \multicolumn{1}{m{3.8cm}|}{$\begin{bmatrix} H_1^{-1} & 0 \\ 0 & H_2^{-1} \end{bmatrix}^{\w{-}} \hspace{-0.2cm}
    \begin{bmatrix} I & 0 \\ -S_{x_1}^{x_2} & I \end{bmatrix}$} \\[0.55cm] \hline
    
    \multicolumn{1}{|c|}{\begin{tabular}{c}
         Block diagram \\ of M in Fig.~\ref{fig:intersysdiagM}
    \end{tabular}} 
    & {
    \begin{tikzpicture}[auto, node distance=2cm,>=latex',baseline=(dx1)]
    \node [input,name=f1] (f1) {};
    \node [tmp,name=tf1,right of = f1, node distance =1cm] (tf1) {};
    \node [block, name=m1, below of = tf1, node distance = 0.5cm] (m1) {$1$};
    \node [block, name=int1, below of = m1, node distance = 1.5cm] (int1) {$\int$};
    \node [tmp,name=tx1,below of = int1, node distance =0.5cm] (tx1) {};
    \node [input,name=x1, left of =tx1, node distance = 1cm] (x1) {};
    
    \draw[->] (f1) |- node[above=-2.5pt,pos=0.6]{$f_1$} (tf1) -| (m1);
    \draw[->] (m1) -- node[left=-2.5pt,pos=0.5,name=dx1]{$\dot{x}_1$}  (int1);
    \draw[->] (int1) |- (tx1) -- node[below=-1.5pt,pos=0.9]{$x_1$} (x1);
    
    \node [block,name=0f2,below of =int1, node distance = 1cm] (0f2) {$x_2^s(x_1)$};
    \node [input,name=x2,below of = x1, node distance = 0.5cm] (x2) {};
    \draw[->] (int1) -- (0f2);
    \draw[->] (0f2) -- node[below=-1.5pt,pos=0.8]{$x_2$} (x2);
\end{tikzpicture} 
\vspace{-0.0cm} }
    & { 
    \begin{tikzpicture}[auto, node distance=2cm,>=latex',baseline=(dx1)]
    \node [input,name=f1] (f1) {};
    \node [tmp,name=tf1,right of = f1, node distance =0.8cm] (tf1) {};
    \node [block, name=m1, below of = tf1, node distance = 0.5cm] (m1) {$1$};
    \node [block, name=int1, below of = m1, node distance = 1.5cm] (int1) {$\int$};
    \node [tmp,name=tx1,below of = int1, node distance =0.5cm] (tx1) {};
    \node [input,name=x1, left of =tx1, node distance = 0.8cm] (x1) {};
    
    \draw[->] (f1) |- node[above=-2.5pt,pos=0.6]{$f_1$} (tf1) -| (m1);
    \draw[->] (m1) -- node[left=-2.5pt,pos=0.5,name=dx1]{$\dot{x}_1$}  (int1);
    \draw[->] (int1) |- (tx1) -- node[below=-1.5pt,pos=0.7]{$x_1$} (x1);
    
    \node [input,name=f2,right of = x1, node distance = 2.8cm] (f2) {};
    \node [tmp,name=tf2,left of = f2, node distance =0.8cm] (tf2) {};
    \node [block, name=m2, above of = tf2, node distance = 0.5cm] (m2) {$\frac{1}{\epsilon}$};
    \node [block, name=int2, above of = m2, node distance = 1.5cm] (int2) {$\int$};
    \node [tmp,name=tx2,above of = int2, node distance =0.5cm] (tx2) {};
    \node [input,name=x2, right of =tx2, node distance = 0.8cm] (x2) {};
    
    \draw[->] (f2) |- node[below=-2.5pt,pos=0.6]{$f_2$} (tf2) -| (m2);
    \draw[->] (m2) -- node[right=-2.5pt,pos=0.5,name=dx2]{$\dot{x}_2$}  (int2);
    \draw[->] (int2) |- (tx2) -- node[above=-2.5pt,pos=0.7]{$x_2$} (x2);
    
    \node [input,name=x1p,above of = x2, node distance = 0.5cm] (x1p) {};
    \node [tmp, name= tx1p3, right of =tx1, node distance = 0.5cm] (tx1p3) {};
    \node [tmp, name= tx1p4, above of =tx1p3, node distance = 3cm] (tx1p4) {};
    \draw[->] (tx1) -- (tx1p3) --(tx1p4) --node[above=-1.5pt,pos=0.8]{$x_1$} (x1p);
 
    \node [input,name=x2p,below of = x1, node distance = 0.5cm] (x2p) {};
    \node [tmp, name= tx2p3, left of =tx2, node distance = 0.5cm] (tx2p3) {};
    \node [tmp, name= tx2p4, below of =tx2p3, node distance = 3cm] (tx2p4) {};
    \draw[->] (tx2) -- (tx2p3) --(tx2p4) --node[below=-1.5pt,pos=0.8]{$x_2$} (x2p);  
\end{tikzpicture}
\vspace{-0.0cm}}
    & {
    \begin{tikzpicture}[auto, node distance=2cm,>=latex',baseline=(dx1)]
    \node [input,name=f1] (f1) {};
    \node [tmp,name=tf1,right of = f1, node distance =0.5cm] (tf1) {};
    \node [block, name=m1, below of = tf1, node distance = 0.5cm,minimum height=1em] (m1) {$1$};
    \node [block, name=int1, below of = m1, node distance = 1.5cm,minimum height=1em] (int1) {$\int$};
    \node [tmp,name=tx1,below of = int1, node distance =0.5cm] (tx1) {};
    \node [input,name=x1, left of =tx1, node distance = 0.5cm] (x1) {};
    
    \draw[->] (f1) |- node[above=-2.5pt,pos=0.6]{$f_1$} (tf1) -| (m1);
    \draw[->] (m1) -- node[left=-2.5pt,pos=0.5,name=dx1]{$\dot{x}_1$}  (int1);
    \draw[->] (int1) |- (tx1) -- node[below=-1.5pt,pos=0.6]{$x_1$} (x1);
    
    \node [input,name=f2,right of = x1, node distance = 3.2cm] (f2) {};
    \node [tmp,name=tf2,left of = f2, node distance =0.5cm] (tf2) {};
    \node [block, name=m2, above of = tf2, node distance = 0.5cm,minimum height=1em] (m2) {$1$};
    \node [block, name=int2, above of = m2, node distance = 1.5cm,minimum height=1em] (int2) {$\int$};
    \node [tmp,name=tx2,above of = int2, node distance =0.5cm] (tx2) {};
    \node [input,name=x2, right of =tx2, node distance = 0.5cm] (x2) {};
    
    \draw[->] (f2) |- node[below=-2.5pt,pos=0.6]{$f_2$} (tf2) -| (m2);
    \draw[->] (int2) |- (tx2) -- node[above=-1.5pt,pos=0.6]{$x_2$} (x2);
    
    \node [input,name=x1p,above of = x2, node distance = 0.5cm] (x1p) {};
    \node [tmp, name= tx1p3, right of =tx1, node distance = 0.5cm] (tx1p3) {};
    \node [tmp, name= tx1p4, above of =tx1p3, node distance = 3cm] (tx1p4) {};
    \draw[->] (tx1) -- (tx1p3) --(tx1p4) --node[above=-1.5pt,pos=0.9]{$x_1$} (x1p);
 
    \node [input,name=x2p,below of = x1, node distance = 0.5cm] (x2p) {};
    \node [tmp, name= tx2p3, left of =tx2, node distance = 0.5cm] (tx2p3) {};
    \node [tmp, name= tx2p4, below of =tx2p3, node distance = 3cm] (tx2p4) {};
    \draw[->] (tx2) -- (tx2p3) --(tx2p4) --node[below=-1.5pt,pos=0.9]{$x_2$} (x2p);   
        
    \node [sum,name=sum,above of = m2, node distance = 0.75cm, fill=orange ,fill opacity=.2, draw opacity=1, text opacity = 1] (sum) {};
    \node [tmp, name=tdx1,below of = m1, node distance =0.75cm] (tdx1) {};
    \node [block,name=S,right of = tdx1,node distance = 1.1cm, fill=orange ,fill opacity=.2, draw opacity=1, text opacity = 1] (S) {$S_{x_1}^{x_2}$};    
    \node [tmp, name = tx1p6, right of = tx1p3, node distance = 0.6cm] (tx1p6) {};
    \node [tmp, name = tx2p6, left of = tx2p3, node distance = 0.6cm] (tx2p6) {};
    \draw[->] (m2) -- (sum);
    \draw[->] (sum) -- node[right=-1pt,pos=0.3,name=dx2]{$\dot{x}_2$} (int2);
    \draw[->,orange,densely dotted,thick] (tdx1) -- (S);
    \draw[->,orange,densely dotted,thick] (S) -- (sum);
    \draw[->,orange,densely dotted,thick] (tx1p3) -- node[below=-2.5pt]{$x_1$} (tx1p6) -- (S) ;
    \draw[->,orange,densely dotted,thick] (tx2p3) -- node[above=-2.5pt]{$x_2$} (tx2p6) -- (S) ;
\end{tikzpicture}
\vspace{-0.0cm}}
    & {
    \begin{tikzpicture}[auto, node distance=2cm,>=latex',baseline=(dx1)]
    \node [input,name=f1] (f1) {};
    \node [tmp,name=tf1,right of = f1, node distance =0.8cm] (tf1) {};
    \node [block, name=m1, below of = tf1, node distance = 0.5cm] (m1) {$H_1$};
    \node [block, name=int1, below of = m1, node distance = 1.5cm] (int1) {$\int$};
    \node [tmp,name=tx1,below of = int1, node distance =0.5cm] (tx1) {};
    \node [input,name=x1, left of =tx1, node distance = 0.8cm] (x1) {};
    
    \draw[->] (f1) |- node[above=-2.5pt,pos=0.6]{$f_1$} (tf1) -| (m1);
    \draw[->] (m1) -- node[left=-2.5pt,pos=0.5,name=dx1]{$\dot{x}_1$}  (int1);
    \draw[->] (int1) |- (tx1) -- node[below=-1.5pt,pos=0.7]{$x_1$} (x1);

    \node [input,name=f2,right of = x1, node distance = 3.8cm] (f2) {};
    \node [tmp,name=tf2,left of = f2, node distance =0.8cm] (tf2) {};
    \node [block, name=m2, above of = tf2, node distance = 0.5cm] (m2) {$H_2$};
    \node [block, name=int2, above of = m2, node distance = 1.5cm] (int2) {$\int$};
    \node [tmp,name=tx2,above of = int2, node distance =0.5cm] (tx2) {};
    \node [input,name=x2, right of =tx2, node distance = 0.8cm] (x2) {};
    
    \draw[->] (f2) |- node[below=-2.5pt,pos=0.6]{$f_2$} (tf2) -| (m2);
    \draw[->] (int2) |- (tx2) -- node[above=-1.5pt,pos=0.7]{$x_2$} (x2);
    
    \node [input,name=x1p,above of = x2, node distance = 0.5cm] (x1p) {};
    \node [tmp, name= tx1p3, right of =tx1, node distance = 0.5cm] (tx1p3) {};
    \node [tmp, name= tx1p4, above of =tx1p3, node distance = 3cm] (tx1p4) {};
    \draw[->] (tx1) -- (tx1p3) --(tx1p4) --node[above=-1.5pt,pos=0.9]{$x_1$} (x1p);
 
    \node [input,name=x2p,below of = x1, node distance = 0.5cm] (x2p) {};
    \node [tmp, name= tx2p3, left of =tx2, node distance = 0.5cm] (tx2p3) {};
    \node [tmp, name= tx2p4, below of =tx2p3, node distance = 3cm] (tx2p4) {};
    \draw[->] (tx2) -- (tx2p3) --(tx2p4) --node[below=-1.5pt,pos=0.9]{$x_2$} (x2p);   
    
    \node [tmp, name=tx1p, left of = tx1, node distance = 0.6cm] (tx1p) {};
    \node [tmp, name=tx1pp, left of = m1, node distance = 0.6cm] (tx1pp) {};
    \node [tmp, name=tx2p, right of = tx2, node distance = 0.6cm] (tx2p) {};
    \node [tmp, name=tx2pp, right of = m2, node distance = 0.6cm] (tx2pp) {};
    \node [tmp,name = tx2p5, below of = tx2p3, node distance = 0.5cm] (tx2p5) {};
    \node [tmp,name = tx1p5, above of = tx1p3, node distance = 0.5cm] (tx1p5) {};
    \draw[->,blue,dashed,thick] (tx1p) |- (tx1pp) -- (m1);
    \draw[->,blue,dashed,thick] (tx2p) |- (tx2pp) -- (m2);
    \draw[->,blue,dashed,thick] (tx2p5) -- node[above=-2.5pt,pos=0.2]{$x_2$} (m1);
    \draw[->,blue,dashed,thick] (tx1p5) -- node[below=-2.5pt,pos=0.2]{$x_1$} (m2);
        
    \node [sum,name=sum,above of = m2, node distance = 0.75cm, fill=orange ,fill opacity=.2, draw opacity=1, text opacity = 1] (sum) {};
    \node [tmp, name=tdx1,below of = m1, node distance =0.75cm] (tdx1) {};
    \node [block,name=S,right of = tdx1,node distance = 1.1cm, fill=orange ,fill opacity=.2, draw opacity=1, text opacity = 1] (S) {$S_{x_1}^{x_2}$};    
    \node [tmp, name = tx1p6, right of = tx1p3, node distance = 0.6cm] (tx1p6) {};
    \node [tmp, name = tx2p6, left of = tx2p3, node distance = 0.6cm] (tx2p6) {};
    \draw[->] (m2) --  (sum);
    \draw[->] (sum) -- node[right=-1pt,pos=0.3,name=dx2]{$\dot{x}_2$} (int2);
    \draw[->,orange,densely dotted,thick] (tdx1) -- (S);
    \draw[->,orange,densely dotted,thick] (S) -- (sum);
    \draw[->,orange,densely dotted,thick] (tx1p3) -- node[below=-2.5pt]{$x_1$} (tx1p6) -- (S) ;
    \draw[->,orange,densely dotted,thick] (tx2p3) -- node[above=-2.5pt]{$x_2$} (tx2p6) -- (S) ;
\end{tikzpicture} 
\vspace{-0.0cm}} \\ \hline
  
    
    
    \begin{tabular}{c}
         Jacobian $J$ at \\ $(x_1^s,x_2^s(x_1^s))$
    \end{tabular}  
    & \multicolumn{1}{m{2.5cm}|}{$\begin{matrix} J_1 = \nabla_{x_1} f^r_1 \\ 
   J_2 = \nabla_{x_2} f_2
    \end{matrix}$ }
    & \multicolumn{1}{m{3.3cm}|}{$\begin{bmatrix} \nabla_{x_1} f_1 & \nabla_{x_2} f_1 \\ \frac{1}{\epsilon}\nabla_{x_1} f_2 & \frac{1}{\epsilon}\nabla_{x_2} f_2
    \end{bmatrix}$ }
    & \multicolumn{1}{m{3.3cm}|}{$\sim \hspace{-0.1cm} \begin{bmatrix} \nabla_{x_1} f^r_1 & \nabla_{x_2}f_1 \\ 0 & \nabla_{x_2}f_2 \end{bmatrix}$} 
    & \multicolumn{1}{m{4cm}|}{$\sim \hspace{-0.1cm} \begin{bmatrix} H_1\nabla_{x_1} f^r_1 & H_1\nabla_{x_2}f_1 \\ 0 & H_2\nabla_{x_1}f_2 \end{bmatrix}$}\\[0.5cm] \hline
    
    
    \begin{tabular}{c}
         Eigenvalues $\lambda$ \\[-0.15cm] \& \\[-0.2cm] local stability 
    \end{tabular} 
    & \multicolumn{1}{m{3cm}|}{\cite[Cor.~4.3]{khalil2002nonlinear}: Exp. stable if and only if $\lambda$ of Jacobians $J_i$ have negative real part.}
    & \multicolumn{1}{m{3.6cm}|}{\textbf{Remark~\ref{rem:stablost}}: different $\lambda$, 
    can be unstable even if \eqref{eq:sys2ideal} stable.
    Stable if $\epsilon < \bar{\epsilon}$.}
    & \multicolumn{1}{m{3.8cm}|}{\textbf{Proposition~\ref{prop:sys2stab2} \& Cor.~\ref{cor:sys2corstab}}: similar to block-diagonal $J$, same $\lambda$ as \eqref{eq:sys2ideal}, thus same local stability.}
    & \multicolumn{1}{m{4.4cm}|}{\textbf{Proposition~\ref{prop:sys2stab4} \& Cor.~\ref{cor:sys2predsensacc}}: Preserving stability, $\lambda$ of $J$ can have lower negative real part than \eqref{eq:sys2ideal} and \eqref{eq:sys2predsens}, thus faster convergence.} \\ \hline
    
  \end{tabular}
\end{table*}

\subsection{Sensitivity-conditioning interconnection}\label{subsec:sys2predsens}

Here we present an alternative interconnection in \eqref{eq:intersys}, that can preserve the steady state $(x_1^s,x_2^s(x_1^s))$ of the two time-scale system \eqref{eq:sys2ideal} and its stability, \rr{without the need of a sufficiently large time-scale separation via a a singular parameter $\epsilon$}.
This interconnection uses a sensitivity-conditioning matrix $M=\left[\begin{smallmatrix} I & 0 \\ -S_{x_1}^{x_2} & I \end{smallmatrix}\right]$, graphically presented in Table~\ref{tab:comparM}:    
\begin{equation}\label{eq:sys2predsens}
      \begin{bmatrix} I & 0 \\ -S_{x_1}^{x_2}(x_1,x_2) & I \end{bmatrix} \hspace{-0.2cm}
     \begin{bmatrix} \dot{x}_1  \\ \dot{x}_2 \end{bmatrix} 
     =\begin{bmatrix} f_1(x_1,x_2) \\ f_2(x_1,x_2) \end{bmatrix} \hspace{-0.1cm}.
\end{equation}

\saverio{Instead of accelerating the second subsystem as in \eqref{eq:sys2singpert}, this conditioning matrix $M$ contains an off-diagonal term that changes the dynamics of the second subsystem to $\dot{x}_2  =f_2(x_1,x_2) + S_{x_1}^{x_2}(x_1,x_2) f_1(x_1,x_2)$, i.e., using the control input $u_2=S_{x_1}^{x_2}(x_1,x_2) f_1(x_1,x_2)$ in \eqref{eq:intersysu}.} Intuitively, there are now two components in the vector field of $\dot{x}_2$: $f_2(x_1,x_2)$ drives $x_2$ to the steady state $x_2^s(x_1)$, while the sensitivity-conditioning $S_{x_1}^{x_2}(x_1,x_2) f_1(x_1,x_2)$ can be interpreted as a \textit{feed-forward} term anticipating the change of $x_2^s(x_1)$ due to the dynamics $\dot{x}_1 \neq 0$. 
\saverio{This second term affects the transient behavior only and vanishes at steady state.}

Given Assumption~\ref{ass:sys2sensinv}, local existence and uniqueness \cite[Thm.~3.1]{khalil2002nonlinear} of a solution $x_i(t)$ for \eqref{eq:sys2predsens} are guaranteed if:
\begin{ass}\label{ass:sys2loclip}
The vector field $f_2(x_1,x_2) + S_{x_1}^{x_2}(x_1,x_2)f_1(x_1,x_2)$ is locally Lipschitz continuous.
\end{ass}
For more insight on the benefits of \eqref{eq:sys2predsens}, we advance some results, that specialize the more general Theorem~\ref{thm:sysNpredsensstab} (presented later in Section~\ref{sec:multits}) to the case of two interconnected systems. The first proposition shows that \rr{the singleton $\{x_2^s(x_1)\}$ is a positively invariant set, i.e.,} once $x_2$ hits the steady state $x_2^s(x_1)$, it remains at $x_2^s(x_1)$ even if $\dot{x}_1 \neq 0$:

\begin{prop}[\saverio{Positive invariance}]\label{prop:sys2stab1} 
Consider the dynamics of $x_2$ in \eqref{eq:sys2predsens} initialized at time $t_0$:
\begin{equation}\label{eq:sys2predsens2}
    \dot{x}_2  = f_2(x_1,x_2) + S_{x_1}^{x_2}(x_1,x_2) \dot{x}_1 \text{ s.t. } x_2(t_0)=x_2^s(x_1(t_0)).
\end{equation}
Then, $x_2(t)=x_2^s(x_1(t))$ is the unique solution on the open domain of existence.
\end{prop}

\begin{proof} 
 First, note that $x_2(t)=x_2^s(x_1(t))$ satisfies \eqref{eq:sys2predsens2}:
\begin{equation*}\begin{array}{rl}
     \dot{x}_2(t)  &
     = f_2(x_1(t),x_2(t)) + S_{x_1}^{x_2}(x_1(t),x_2(t)) \dot{x}_1(t) \\
    & = \cancel{f_2(x_1(t),x_2^s(x_1(t)))} + \underbrace{S_{x_1}^{x_2}(x_1(t),x_2^s(x_1(t)))}_{\overset{\eqref{eq:sys2sens}}{=} \nabla_{x_1} x_2^s(x_1(t))} \dot{x}_1(t) \\[-0.4cm]
    & = \frac{dx_2^s(x_1(t))}{dt}
\end{array}
\end{equation*}
Local existence and uniqueness of a solution is guaranteed by Assumption~\ref{ass:sys2loclip}. Thus, $x_2(t)=x_2^s(x_1(t))$ is the unique solution on this domain of existence, because the derivatives and initial conditions of $x_2(t)$ and $x_2^s(x_1(t))$ coincide for $t\geq t_0$.
\end{proof}

Moreover, the sensitivity-conditioning \eqref{eq:sys2predsens} allows to preserve the local stability of the two-time-scale system \eqref{eq:sys2ideal}:

\begin{prop}[Local stability]\label{prop:sys2stab2} 
At a steady state $(x_1^s,x_2^s(x_1^s))$, the Jacobian $J$ of \eqref{eq:sys2predsens} satisfies:
\begin{equation}\label{eq:sys2predsensjac}
     J \sim \begin{bmatrix} \nabla_{x_1}f^r_1(x_1^s) & \nabla_{x_2}f_1(x_1^s,x_2^s(x_1^s)) \\ 0 & \nabla_{x_2}f_2(x_1^s,x_2^s(x_1^s)) \end{bmatrix} \hspace{-0.1cm},
\end{equation}
where $f_1^r$ is the reduced vector field from \eqref{seq:sys2idealro}, and $\sim$ denotes similarity, i.e., related by a similarity transformation. 
\end{prop}

\begin{proof}
To calculate the Jacobian $J$ of \eqref{eq:sys2predsens} at $(x_1^s,x_2^s(x_1^s))$, we invert $M$ as $M^{-1}=\left[\begin{smallmatrix} I & 0 \\ -S_{x_1}^{x_2} & I \end{smallmatrix}\right]^{-1}=\left[\begin{smallmatrix} I & 0 \\ S_{x_1}^{x_2} & I \end{smallmatrix}\right]$, take derivatives, and evaluate them at steady-state, so that $f_i(x_1^s,x_2^s(x_1^s))=0$:
\begin{equation*}\begin{array}{rl}
    J = & M^{-1}
    \begin{bmatrix} \nabla_{x_1}f_1 & \nabla_{x_2}f_1 \\ \nabla_{x_1}f_2 & \nabla_{x_2}f_2 \end{bmatrix} \hspace{-0.1cm},
\end{array}
\end{equation*}
where for clarity we omit the evaluation point $(x_1^s,x_2^s(x_1^s))$ in the notation. This $J$ is similar to $\tilde{J}:= M JM^{-1}$, where
\begin{equation*}
     J \hspace{-0.05cm} \sim \hspace{-0.05cm} \tilde{J} =
     \cancel{ M M^{-1}} \begin{bmatrix} \nabla_{x_1}f_1 & \nabla_{x_2}f_1 \\ \nabla_{x_1}f_2 & \nabla_{x_2}f_2 \end{bmatrix} 
     M^{-1}  \hspace{-0.05cm} \overset{\eqref{eq:sys2sens},\eqref{eq:sys2totder}}{=} \hspace{-0.05cm} \begin{bmatrix} \nabla_{x_1}f^r_1 & \nabla_{x_2}f_1 \\ 0 & \nabla_{x_2}f_2 \end{bmatrix} \hspace{-0.1cm}.
\end{equation*}
\end{proof}

\begin{cor}\label{cor:sys2corstab}
The Jacobians of the two time-scale system \eqref{eq:sys2ideal} and the sensitivity-conditioning conditioned system \eqref{eq:sys2predsens} have the same eigenvalues, and thus the same local stability properties. 
\end{cor}

\rr{
\begin{rem}\label{rem:crosscanc}
The cancellation of one off-diagonal term in \eqref{eq:sys2predsensjac} is due to the sensitivity definition in \eqref{eq:sys2sens}, and will also play a crucial role in the proofs of the results to come. Essentially, the role of the sensitivity-conditioning is to turn a closed-loop into a cascade system from the viewpoint of stability analysis, see also the later Remark~\ref{rem:backstep}.
\end{rem}
}


\rr{
Propositions~\ref{prop:sys2stab1} and \ref{prop:sys2stab2} and Corollary~\ref{cor:sys2corstab} establish that invariance and local exponential stability of $x_2^s(x_1)$ are preserved from the two-time-scale system \eqref{eq:sys2ideal} in the single-time-scale sensitivity-conditioning one \eqref{eq:sys2predsens}. Furthermore, these results can be extended to contraction regions satisfying the following:
\begin{ass}[Contraction region \cite{lohmiller1998contraction}]\label{ass:contreg}
There exists $\eta_2>0$ and an open ball $\mathcal{B}_{r_2}\big(x_2^s(x_1)\big)=\{x_2 | \normsz{x_2-x_2^s(x_1)}_2 < r_2\}$ centered at $x_2^s(x_1)$, with a positive radius $r_2>0$, and a metric defined by a constant symmetric positive definite $P_2 \succ 0$, such that for $x_2 \in \mathcal{B}_{r_2}\big(x_2^s(x_1)\big)$ it holds uniformly for all $x_1$ that
$$P_2\nabla_{x_2}f_2(x_1,x_2) + \nabla_{x_2}f_2(x_1,x_2)^T P_2 \preceq -\eta_2 P_2$$ 
\end{ass}
}

\rr{Under this Assumption~\ref{ass:contreg}, $\mathcal{B}_{r_2}\big(x_2^s(x_1)\big)$ is a contraction region for the boundary-layer system \eqref{seq:sys2ideal0bl} for all $x_1$, within which the invariant set $x_2^s(x_1)$ (see Proposition~\ref{prop:sys2stab1}) is exponentially stable \cite[Thm.~2]{lohmiller1998contraction}. Then, such a contraction region is also preserved under the sensitivity-conditioning \eqref{eq:sys2predsens}:}

\rr{
\begin{prop}[Stability with a contracting boundary layer]\label{prop:sys2stab3}
Under Assumption~\ref{ass:contreg} the following holds:
\begin{enumerate}[leftmargin=*] 
    \item\label{itempf1} The sensitivity-conditioning interconnection \eqref{eq:sys2predsens} is well-defined for $x_2 \in \mathcal{B}_{r_2}(x_2^s(x_1))$, i.e., $\nabla_{x_2}f_2(x_1,x_2)$ is uniformly invertible in $\mathcal{B}_{r_2}(x_2^s(x_1))$. 
    \newline
    Furthermore, there exists $r_{2,0}<r_2$ such that if $x_2(0) \in \mathcal{B}_{r_{2,0}}\big(x_2^s(x_1(0))\big)$, then $x_2(t) \in \mathcal{B}_{r_2}\big(x_2^s(x_1(t))\big)$ for all $t>0$, and $x_2^s(x_1(t))$ is a locally exponentially stable trajectory for $x_2$ in \eqref{eq:sys2predsens}, i.e., there exists $\eta,K>0$ so that $\normsz{x_2(t)-x_2^s(x_1(t))}_2 \leq K\normsz{x_2(0)-x_2^s(x_1(0))}_2e^{-\eta t}$. 
    \item\label{itempf2} Additionally, assume that $x_1^s$ is an asymptotically stable steady state of the reduced-order system \eqref{seq:sys2idealro}, that the ball $\mathcal{B}_{r_1}(x_1^s)$ is in its region of attraction, and that $f_1^r(\cdot)$ is continuously differentiable in the closure of $\mathcal{B}_{r_1}(x_1^s)$. 
    \newline
    Then, there exists $r_{1,0}\leq r_1$ and $\tilde{r}_{2,0} \leq {r}_{2,0}$, such that if $x_1(0) \in \mathcal{B}_{{r}_{1,0}}(x_1^s)$ and $x_2(0) \in \mathcal{B}_{\tilde{r}_{2,0}}(x_2^s(x_1(0)))$, then $x_1(t) \in \mathcal{B}_{r_1}(x_1^s)$, and $x_1^s$ is asymptotically stable under the sensitivity-conditioning interconnection \eqref{eq:sys2predsens}. 
\end{enumerate} 
\end{prop}
}

\begin{proof} 
\rr{We use the following technical result:
\begin{lemma}\label{lem:propext}
Consider a system $\dot{x}=f(x)$, with steady-state $x^s$ and a continuous differentiable $f(\cdot)$. If there exist a radius $r>0$, a symmetric positive definite matrix $P\succ 0$, and a parameter $\eta>0$, such that $P\nabla_{x}f(x) + (\nabla_{x}f(x))^T P^T \preceq -\eta P$ for $x \in \mathcal{B}_{r}(x^s)$, then in $\mathcal{B}_{r}(x^s)$ it holds that
\begin{enumerate}[leftmargin=*]
    \item the inverse of the Jacobian $\nabla_{x}f(x)$ exists and is bounded: $\normsz{\nabla_{x}f(x)^{-1}}_2 \leq \tfrac{2\lambda_{\max}(P) }{\eta\lambda_{\min}(P)}$, and
    \item the vector field $f(x)$ is lower bounded: \newline $ \normsz{f(x)}_P \geq \tfrac{\eta}{2} \normsz{x-x^s}_P$.
\end{enumerate}
Here $\lambda_{\min}(\cdot), \allowbreak \lambda_{\max}(\cdot)$ denote the minimum and maximum eigenvalues.
\end{lemma}
\begin{proof}
See Appendix~\ref{app:propextproof}. We remark that the second result can be seen as a particular case of \cite[Prop.~3]{aminzare2014contraction}.
\end{proof}
\ref{itempf1}) First, note that Lemma~\ref{lem:propext} assures non-singularity of $\nabla_{x_2}f_2$. Consider the following Krasovskii Lyapunov function $V_2(x_1,x_2)=\normsz{f_2(x_1,x_2)}_{P_2}^2$ \cite[Ch.~5]{sastry2013nonlinear} for $x_2$ in \eqref{eq:sys2predsens}. Since $P_2 \nabla_{x_2} f_2 + \nabla_{x_2}f_2^T P_2^T \preceq -\eta_2 P_2$, under the sensitivity-conditioning dynamics \eqref{eq:sys2predsens} we have:
\begin{equation}\label{eq:V2dot}\begin{array}{rcl}     
\dot{V}_2 
    & = & f_2^T P_2\nabla_{x_1}f_2 \dot{x}_1 + f_2^T P_2 \nabla_{x_2}f_2 \dot{x}_2 \\
    & &  + (f_2^T P_2 \nabla_{x_1}f_2 \dot{x}_1 + f_2^T  P_2 \nabla_{x_2}f_2 \dot{x}_2)^T \\
    & \overset{\eqref{eq:sys2predsens}}{=} & f_2^T P_2 \big(\cancel{\nabla_{x_1}f_2 +  \nabla_{x_2}f_2 S_{x_1}^{x_2}} \big) \dot{x}_1 + f_2^T P_2 \nabla_{x_2}f_2 f_2 \\
    & &  + (f_2^T P_2 \big(\cancel{\nabla_{x_1}f_2 +  \nabla_{x_2}f_2 S_{x_1}^{x_2}} \big) \dot{x}_1 + f_2^T P_2 \nabla_{x_2}f_2 f_2 )^T \\
    & \overset{\eqref{eq:sys2sens}}{=} & f_2^T \big(P_2 \nabla_{x_2} f_2 + \nabla_{x_2}f_2^T P_2^T \big) f_2 \leq -\eta_2 V_2 
\end{array}
\end{equation}
Hence, $\normsz{f_2(x_1(t),x_2(t))}_{P_2}^2 \leq \normsz{f_2(x_1(0),x_2(0))}_{P_2}^2 e^{-\eta_2 t}$. Since $f_2(\cdot)$ is continuously differentiable, it is locally Lipschitz continuous in $\mathcal{B}_{r_2}(x_2^s(x_1))$ with some constant $L_{f_2}$, and by involving Lemma~\ref{lem:propext} we have
\begin{equation}\label{eq:convx2}\begin{array}{l}
     \normsz{x_2(t)-x_2^s(x_1(t))}_{2}^2 \leq \frac{4}{\eta_2^2\lambda_{\min}(P_2)} \normsz{f_2(x_2(t),x_1(t))}_{P_2}^2  \\
\leq \underbrace{\tfrac{4L_{f_2}^2\lambda_{\max}(P_2)}{\eta_2^2\lambda_{\min}(P_2)}}_{=\tfrac{r_2^2}{r_{2,0}^2}}e^{-\eta_2t} \normsz{x_2(0)  - x_2^s(x_1(0))}_{2}^2,    
\vspace{-0.5cm}
\end{array}
\end{equation}
where $r_{2,0}:=\tfrac{\eta_2r_2}{2 L_{f_2}}\sqrt{\tfrac{\lambda_{\min}(P_2)}{\lambda_{\max}(P_2)}} \leq r_2$. Hence, if $\normsz{x_2(0)-x_2^s(x_1(0))}_{2} \allowbreak < r_{2,0}$, then $x_2(t) \in \mathcal{B}_{r_2}\big(x_2^s(x_1(t))\big)$ for all $t>0$, and $x_2$ converges exponentially to $x_2^s(x_1)$ under \eqref{eq:sys2predsens}, despite the varying $x_1$.
}

\rr{
\ref{itempf2}) Now we analyse the $x_1$-dynamics subject to the exponential converging input $x_2(t)-x_2^s(x_1(t))$. Since $f_1^r(\cdot)$ is continuously differentiable, it is locally Lipschitz continuous, and $\nabla_{x_1}f_1^r(\cdot)$ is bounded in $\mathcal{B}_{r_1}(x_1^s)$. Hence, by the converse Lyapunov theorem \cite[Thm.~4.16]{khalil2002nonlinear}, there exists a Lyapunov function $V_1(x_1)$ satisfying:
\begin{equation}\label{eq:lyapro}\begin{array}{c}
     \alpha_1(\normsz{x_1-x_1^s}_2) \leq V_1(x_1) \leq \alpha_2(\normsz{x_1-x_1^s}_2)  \\
     \nabla_x V_1(x_1)^T f_1^r(x_1) \leq -\alpha_3(\normsz{x_1-x_1^s}_2) \\
     \normsz{\nabla_x V_1(x_1)}_2 \leq \alpha_4(\normsz{x_1-x_1^s}_2),
\end{array}
\end{equation}
where $\alpha_i(\cdot)$ are $\mathcal{K}$-functions. Since $f_1(\cdot)$ is continuously differentiable, it is locally Lipschitz continuous in $\mathcal{B}_{r_1}(x_1^s)\times \mathcal{B}_{r_2}(x_2^s(x_1))$ with some constant $L_{f_1}$. Then, the Lyapunov function $V_1(x_1)$ under the sensitivity-conditioning dynamics \eqref{eq:sys2predsens} satisfies:
\begin{equation*}\begin{array}{rcl}
    \dot{V}_1 & \leq & \nabla_{x_1} V_1^T f_1 \leq \nabla_{x_1} V_1^T f_1^r + \normsz{\nabla_x V_1(x_1)}_2\normsz{f_1-f_1^r}_2 \\
    & \stackrel{\eqref{eq:lyapro}}{\leq} & \hspace{-0.1cm} -\alpha_3(\normsz{x_1 \hspace{-0.05cm} - \hspace{-0.05cm} x_1^s}_2)  + L_{f_1}\alpha_4(\normsz{x_1 \hspace{-0.05cm} - \hspace{-0.05cm} x_1^s}_2)\normsz{x_2 \hspace{-0.05cm} - \hspace{-0.05cm} x_2^s(x_1)}_2 \\
    & \stackrel{\eqref{eq:convx2}}{\leq} & \hspace{-0.1cm} -\alpha_3(\alpha_2^{-1}(V_1)) \\
    & & + \alpha_4(\alpha_1^{-1}(V_1))\frac{r_2 L_{f_1}}{r_{2,0}}e^{-\eta_2 t}\normsz{x_2(0)-x_2^s(x_1(0))}_{2}.
\end{array}
\end{equation*}
Consider any $\delta \in (0,r_1)$, and define $r_{1,0} := \alpha_2^{-1}(\alpha_1(r_1-\delta)) < r_1$ and $\tilde{r}_{2,0} :=  \min \Big( r_{2,0}, \allowbreak \frac{\alpha_3(r_{1,0})}{\alpha_4(r_{1}-\delta)\frac{r_2 L_{f_1}}{r_{2,0}}} \Big)$. If $x_1(0) \in \mathcal{B}_{{r}_{1,0}}(x_1^s)$, then $V_1(x_1(0)) \leq \alpha_1(r_1-\delta)$; and if $x_2(0) \in \mathcal{B}_{\tilde{r}_{2,0}}\big(x_2^s(x_1(0))\big)$, then $\dot{V}_1 \leq 0$ whenever $V_1 = \alpha_1(r_1-\delta)$. Hence, $V_1(x_1(t)) < \alpha_1(r_1)$ and $\normsz{x_1(t)-x_1^s}_2 < r_1$ for all $t>0$. Furthermore, for any $\epsilon \in (0,1)$ it holds that
\begin{equation*}\begin{array}{rcl}
     \dot{V}_1 & \leq& \hspace{-0.1cm} -\alpha_3(\normsz{x_1 \hspace{-0.05cm} - \hspace{-0.05cm} x_1^s}_2)  + L_{f_1}\alpha_4(\normsz{x_1 \hspace{-0.05cm} - \hspace{-0.05cm} x_1^s}_2)\normsz{x_2 \hspace{-0.05cm} - \hspace{-0.05cm} x_2^s(x_1)}_2 \\
     & \leq & \hspace{-0.1cm} -\epsilon\alpha_3(\normsz{x_1 \hspace{-0.05cm} - \hspace{-0.05cm} x_1^s}_2) - (1-\epsilon)\alpha_3(\normsz{x_1 \hspace{-0.05cm} - \hspace{-0.05cm} x_1^s}_2)\\
     & &   + L_{f_1}\alpha_4(r_1)\frac{r_2 L_{f_1}}{r_{2,0}}e^{-\eta_2 t}\normsz{x_2(0)-x_2^s(x_1(0))}_{2} \\
     & \leq & - \epsilon \alpha_3(\normsz{x_1-x_1^s}_2),
\end{array}
\end{equation*}
where the last inequality holds
while $\normsz{x_1-x_1^s}_2 \geq \alpha_3^{-1} \big( \tfrac{1}{1-\epsilon} L_{f_1}\alpha_4(r_1)\frac{r_2 L_{f_1}}{r_{2,0}}e^{-\eta_2 t}\normsz{x_2(0)-x_2^s(x_1(0))}_{2} \big)$. Hence, 
$x_1^s$ is asymptotically stable, because it is input-to-state stable \cite[Thm.~4.18]{khalil2002nonlinear} with respect to a vanishing input, see also \cite[Lemma~4.7]{khalil2002nonlinear}. 
}
\end{proof}

\rr{
\begin{rem}[Connection to contraction theory]\label{rem:contraction}
The exponential stability of the boundary-layer system \eqref{seq:sys2ideal0bl} for a constant $x_1$ is a standard assumption in the context of singular perturbation analysis \cite[Thm.~11.4 and after]{khalil2002nonlinear}, and exponential stability implies the existence of a contraction region \cite[Reverse Thm.~2]{lohmiller1998contraction}. Proposition~\ref{prop:sys2stab3} establishes that the boundary-layer exponential stability can be preserved in the single-time-scale interconnection \eqref{eq:intersys} using the sensitivity-conditioning \eqref{eq:sys2predsens}, independently of $\dot{x}_1$. Essentially, the sensitivity-conditioning \eqref{eq:sys2predsens} turns a system $\Sigma_2$ that is only contracting under a constant $x_1$, as in \eqref{seq:sys2ideal0bl}, into a \textit{partially contracting} system \eqref{eq:sys2predsens} in $x_2$ \cite[Def.~1]{del2012contraction} under a time-varying $x_1$. Then, $\normsz{x_2-x_2^s(x_1)}$ becomes an exponentially decaying perturbation for $x_1$ in \eqref{eq:sys2predsens}, and thus asymptotic stability of the reduced-order system \eqref{seq:sys2idealro} can be preserved in \eqref{eq:sys2predsens} under some additional conditions.
\end{rem}
\begin{rem}[Connection to backstepping]\label{rem:backstep}
The role of the sensitivity term $S_{x_1}^{x_2}$ in \eqref{eq:sys2predsens} is to cancel a cross term in the stability analysis of $x_2$ that appears under a time-varying $x_1$, 
see \eqref{eq:V2dot} and the proofs of Propositions~\ref{prop:sys2stab2} and \ref{prop:sys2stab3}. In other words, the sensitivity-conditioning \eqref{eq:sys2predsens} is turning an interconnected system \eqref{eq:intersys} into a cascaded one from the viewpoint of stability analysis, see Remark~\ref{rem:crosscanc}. In a more general setting for Proposition~\ref{prop:sys2stab3}, we could assume that the boundary-layer system \eqref{seq:sys2ideal0bl} is asymptotically stable with a general Lyapunov function $V_2(x_1,x_2)$ that is positive definite with respect to $\normsz{x_2-x_2^s(x_1)}$, and satisfies $\nabla_{x_2}V_2^T f_2 \leq 0$ uniformly over $x_1$. Then, under the sensitivity-conditioning dynamics \eqref{eq:sys2predsens} we would have $\dot{V}_2 \leq \nabla_{x_2}V_2^T f_2 + (\nabla_{x_2}V_2^T S_{x_2}^{x_1} - \nabla_{x_1}V_2^T)\dot{x}_1 \leq (\nabla_{x_2}V_2^T S_{x_1}^{x_2} - \nabla_{x_1}V_2^T)\dot{x}_1$, which can be cancelled by choosing an appropriate sensitivity $S_{x_2}^{x_1}$. In Proposition~\ref{prop:sys2stab3} we considered $V_2=\normsz{f_2}_{P_2}^2$. An other option would be $V_2=\normsz{x_2 - x_2^s(x_1)}_{P_2}^2$, which requires the alternative sensitivity $\nabla_{x_1} x_2^s(x_1) \hspace{-0.05cm} = \hspace{-0.05cm} - \hspace{-0.05cm} \nabla_{x_2}f_2(x_1,x_2^s(x_1)\hspace{-0.05cm})^{-1}  \nabla_{x_1}f_2(x_1,x_2^s(x_1)\hspace{-0.05cm})$ to cancel the term $\nabla_{x_2}V_2^T S_{x_2}^{x_2} - \nabla_{x_1}V_2^T = 2(x_2 - x_2^s(x_1)){P_2}(S_{x_2}^{x_2}-\nabla_{x_1}{x_2^s(x_1)})$. This sensitivity $\nabla_{x_1}{x_2^s(x_1)}$ can be interpreted as a backstepping-like approach \cite[Ch.~14]{khalil2002nonlinear} to cancel the dynamics of $\dot{x}_1 \neq 0$ in $\dot{V}_2$, and again turn \eqref{eq:intersys} in a cascaded system from the viewpoint of stability analysis. Note that for this sensitivity similar results as in Propositions~\ref{prop:sys2stab2} and \ref{prop:sys2stab3} can be derived under suitable assumptions. However, implementing the corresponding interconnection $M=\left[\begin{smallmatrix} I & 0 \\ -\nabla_{x_1}{x_2^s(x_1)} & I \end{smallmatrix}\right]$ (or equivalently $u_2 = \nabla_{x_1}{x_2^s(x_1)} f_1(x_1,x_2) = \frac{d}{dt}x_2^s(x_1(t))$) may not be feasible, since it requires a closed-form expression for $x_2^s(x_1)$ to evaluate $\nabla_{x_1}{x_2^s(x_1)}$. Such a closed-form expression for $x_2^s(x_1)$ may be available in special cases, see cascade control in Section~\ref{sec:casccont}, but not in general, see bilevel optimization in Section~\ref{sec:bilopt}. On the other hand, choosing $V_2=\normsz{f_2}_{P_2}^2$ results in the sensitivity in \eqref{eq:sys2sens}, which does not require to know $x_2^s(x_1)$. In this context, the sensitivity-conditioning \eqref{eq:sys2predsens} acts as an implementable substitute for such a backstepping-like approach, with the same local properties. 
\end{rem}
}


\rr{
\begin{cor}[Global exponential stability]\label{cor:sys2stab3}
Assume that the vector fields $f_i(x_1,x_2)$ are Lipschitz continuous, and there exist $P_i\succ 0$ and $\eta_i>0$ such that the following contraction conditions hold globally for all $x_1$ and $x_2$:
\begin{equation}\label{eq:sys2condstab}
    \begin{array}{rl}
     {P_1\nabla_{x_1}f_1^r(x_1) + (\nabla_{x_1}f_1^r(x_1))^T P_1} & \preceq -\eta_1 P_1  \\
     {P_2\nabla_{x_2}f_2(x_1,x_2)+\nabla_{x_2}f_2(x_1,x_2)^T P_2} & \preceq -\eta_2 P_2
\end{array}
\end{equation}
Then $(x_1^s,x_2^s(x_1^s))$ is a globally exponentially stable steady state of both the two time-scale system \eqref{eq:sys2ideal} and the sensitivity-conditioning interconnection \eqref{eq:sys2predsens}.
\end{cor}
\begin{proof} 
If the contraction conditions \eqref{eq:sys2condstab} hold, \eqref{eq:sys2ideal} is globally exponentially stable \cite[Thm.~2]{lohmiller1998contraction}. 
For \eqref{eq:sys2predsens}, consider now the Lyapunov function $V(x_1,x_2)= V_1(x_1) + \theta V_2(x_1,x_2)$, where $\theta>0$, $V_1(x_1)=\normsz{f_1^r(x_1)}_{P_1}^2$ and $V_2(x_1,x_2)=\normsz{f_2(x_1,x_2)}_{P_2}^2$. From \eqref{eq:V2dot} we have $\dot{V}_2 \leq -\eta_2 V_2$, then
\begin{equation*}\begin{array}{rcl}     
\dot{V_1} 
     & = & 
(f_1^r)^T P_1 \nabla_{x_1}f_1^r \dot{x}_1 + ((f_1^r)^T P_1 \nabla_{x_1}f_1^r \dot{x}_1)^T \\
   & \overset{\eqref{eq:sys2predsens}}{=} & 
   (f_1^r)^T \big( P_1 \nabla_{x_1}f_1^r + (\nabla_{x_1}f_1^r)^T P_1^T \big) f_1^r \\
   & & 
    + 2(f_1^r)^T  P_1 \nabla_{x_1}f_1^r (f_1  -f_1^r) \\ 
    & \overset{\eqref{eq:sys2condstab}}{\leq} & 
    - \eta_1 \normsz{f_1^r}_{P_1}^2 \hspace{-0.1cm} + 2L_{f_1} \normsz{f_1^r}_{P_1} \normsz{f_1-f_1^r}_{P_1} \\
    & \leq & 
    - \eta_1 \normsz{f_1^r}_{P_1}^2 \hspace{-0.1cm} + 2L_{f_1}^2 \sqrt{\tfrac{\lambda_{\max}(P_1)}{\lambda_{\min}(P_2)}} \normsz{f_1^r}_{P_1} \normsz{x_2-x_2^s}_{P_2} \\
    & \stackrel{Lem.~\ref{lem:propext}}{\leq} & 
    - \eta_1 \normsz{f_1^r}_{P_1}^2 \hspace{-0.1cm} + \underbrace{ \sqrt{\tfrac{\lambda_{\max}(P_1)}{\lambda_{\min}(P_2)}} \tfrac{4L_{f_1}^2}{\eta} }_{2\nu} \normsz{f_1^r}_{P_1} \normsz{f_2}_{P_2}, 
    \vspace{-0.3cm}
\end{array}
\end{equation*}
Hence,
$\begin{array}{l}
     \dot{V} \stackrel{\eqref{eq:V2dot}}{\leq} - \begin{bmatrix} \normsz{f_2}_{P_2} \\ \normsz{f_1^r}_{P_1} \end{bmatrix}^T
   \hspace{-0.05cm}  \begin{bmatrix} \theta \eta_2 & -\nu \\ -\nu & \eta_1 \end{bmatrix} \hspace{-0.05cm} 
    \begin{bmatrix} \normsz{f_2}_{P_2} \\ \normsz{f_1^r}_{P_1} \end{bmatrix} \hspace{-0.1cm} \leq -\zeta V,
\end{array}
$
for some $\zeta>0$, by choosing $\theta>\frac{\nu^2 }{\eta_1\eta_2}$.
\end{proof}
}

\subsection{Accelerated sensitivity-conditioning}

\rr{The design of the conditioning matrix $M(x_1,x_2)$ in \eqref{eq:intersys} offers to generalize the sensitivity-conditioning \eqref{eq:sys2predsens} to introduce additional degrees of freedom and achieve a better performance of the interconnection, e.g., a faster convergence.}
Consider two \rr{uniformly positive definite matrices $H_1(x_1,x_2)\succ 0,H_2(x_1,x_2) \succ 0$}, and a generalized sensitivity-conditioning:
\begin{equation}\label{eq:sys2predsensacc}\begin{array}{c}
H(x_1,x_2)^{-1}  \hspace{-0.1cm}
      \begin{bmatrix} I & 0 \\ -S_{x_1}^{x_2}(x_1,x_2) & I \end{bmatrix} \hspace{-0.2cm}
     \begin{bmatrix} \dot{x}_1  \\ \dot{x}_2 \end{bmatrix} 
          \hspace{-0.1cm} = \hspace{-0.1cm} \begin{bmatrix} f_1(x_1,x_2) \\ f_2(x_1,x_2) \end{bmatrix} \hspace{-0.1cm}, \\
          
\end{array}
\end{equation}
where $H(x_1,x_2) = \left[\begin{smallmatrix} H_1(x_1,x_2) & 0 \\ 0 & \hspace{-0.3cm} H_2(x_1,x_2) \end{smallmatrix}\right]$.

\begin{prop}[Extension of Propositions~\ref{prop:sys2stab1}, \ref{prop:sys2stab2}, and Corollary~\ref{cor:sys2stab3}]\label{prop:sys2stab4}
The generalized sensitivity-conditioning \eqref{eq:sys2predsensacc} satisfies:
\begin{enumerate}[leftmargin=*]
    \item \rr{The singleton $\{x_2^s(x_1)\}$ is a positively invariant set under the sensitivity-conditioning dynamics \eqref{eq:sys2predsensacc}.}
    \item At steady state $(x_1^s,x_2^s(x_1^s))$, the Jacobian of \eqref{eq:sys2predsensacc} satisfies 
        \begin{equation*}
        J \sim \begin{bmatrix} H_1\nabla_{x_1}f^r_1 & H_1\nabla_{x_2}f_1 \\ 0 & H_2\nabla_{x_2}f_2 \end{bmatrix} \hspace{-0.1cm}
        \end{equation*}
    \item If the vector fields $f_i(x_1,x_2)$ are Lipschitz continuous, and there exists $P_i\succ 0$ and $\eta_i>0$ such that: 
\begin{equation*}
    \begin{array}{rl}
     {P_1 H_1 \nabla_{x_1}f_1^r + (\nabla_{x_1}f_1^r)^T H_1^T P_1} & \preceq -\eta_1 P_1  \\
     {P_2H_2\nabla_{x_2}f_2+\nabla_{x_2}f_2^T H_2^T P_2} & \preceq - \eta_2 P_2, 
\end{array}
\end{equation*}
then $(x_1^s,x_2^s(x_1^s))$ is a globally exponentially stable steady state of the generalized sensitivity-conditioning \eqref{eq:sys2predsensacc}.
\end{enumerate}
For clarity we omit the evaluations at $(x_1^s,x_2^s(x_1^s))$. 
\end{prop}
\begin{proof}
The proof follows analogous steps as the ones for Proposition~\ref{prop:sys2stab1}, \ref{prop:sys2stab2}, and Corollary~\ref{cor:sys2stab3}.
\end{proof}

\begin{cor}[Accelerated sensitivity-conditioning]\label{cor:sys2predsensacc}
If the two-time-scale system \eqref{eq:sys2ideal} is locally exponentially stable, i.e., $J_1=\nabla_{x_1}f_1^r(x_1^s)$ and $J_2=\nabla_{x_2}f_2(x_1^s,x_2^s(x_1^s))$ have eigenvalues with strictly negative real part, the generalized sensitivity-conditioning system \eqref{eq:sys2predsensacc} is locally exponentially stable if using positive scalars $h_i>0$ and $H_i(x_1,x_2)=h_i I$. Moreover, the exponential convergence rate is improved for $h_i>1$. 
\end{cor}


See Table~\ref{tab:comparM} for a comparison of the sensitivity-conditioning approach \eqref{eq:sys2predsens} and \eqref{eq:sys2predsensacc}, the two-time-scale system \eqref{eq:sys2ideal} and the singular pertubated \eqref{eq:sys2singpert}, summarizing these results. 

To conclude, the sensitivity-conditioning \eqref{eq:sys2predsens} allows to preserve the stability of the two-time-scale system \eqref{eq:sys2ideal} in a single time-scale.
This way, the need of artificially slowing down one subsystem (and, consequently, their interconnection) through a singular perturbation \eqref{eq:sys2singpert} is removed. However, a disadvantage of the sensitivity-conditioning is that it 
could produce large inputs $u$ for $\Sigma_2$ in \eqref{eq:intersysu}, even larger if using the generalization in \eqref{eq:sys2predsensacc}, which changes the sensitivity-conditioning term to $S_{x_1}^{x_2}(x_1,x_2) H_1(x_1,x_2)f_1(x_1,x_2)$. This sensitivity-conditioning term could even become unrealizable in systems with control saturation in $u$. On the other hand, if the two-time-scale system \eqref{eq:sys2ideal} is lcoally stable, the generalization \eqref{eq:sys2predsensacc} can be chosen as $H_1(x_1,x_2)={\epsilon}$ with a sufficiently small $\epsilon$, see Corollary~\ref{cor:sys2predsensacc}. Then, $\Sigma_1$ can be slowed down as with a singular perturbation term, and the sensitivity-conditioning term can be made realizable. This interpretation suggests that singular perturbation \eqref{eq:sys2singpert} and sensitivity-conditioning \eqref{eq:sys2predsens} are not mutually exclusive, but can be combined. Interestingly, Corollary~\ref{cor:sys2predsensacc} also allows to choose arbitrary time scales, e.g., $\Sigma_1$ faster than $\Sigma_2$ by choosing $h_1 \gg h_2$, and still preserve the stability of \eqref{eq:sys2ideal}.

\subsection{Robust sensitivity-conditioning}

The sensitivity-conditioning system \eqref{eq:sys2predsens} requires a precise knowledge of the vector fields, essentially the model of the system, to evaluate the matrix $S_{x_1}^{x_2}(x_1,x_2)$. Here we analyse the implications of model errors: assume that instead of $S_{x_1}^{x_2}(x_1,x_2)$, only an approximation $\hat{S}_{x_1}^{x_2}(x_1,x_2)$ is available, and consider the approximated sensitivity-conditioning
\begin{equation}\label{eq:sys2predsensapprox}
      \begin{bmatrix} I & 0 \\ -\hat{S}_{x_1}^{x_2}(x_1,x_2) & I \end{bmatrix} \hspace{-0.2cm}
     \begin{bmatrix} \dot{x}_1  \\ \dot{x}_2 \end{bmatrix} 
     =\begin{bmatrix} f_1(x_1,x_2) \\ f_2(x_1,x_2) \end{bmatrix} \hspace{-0.1cm},
\end{equation}
which also preserves the steady-state $(x_1^s,x_2^s(x_1^s))$.

Under Assumptions~\ref{ass:sys2sensinv} and \ref{ass:sys2loclip}, we can derive a strong robustness certificate in the form of input-to-state stability \cite[Def.~4.7]{khalil2002nonlinear}:
\begin{prop}[Input-to-state stability]\label{prop:sys2rob}
\rr{Consider the same conditions as in Corollary~\ref{cor:sys2stab3}, and define the error
\begin{equation*}\label{eq:sys2predsenserror}
\xi:= (S_{x_1}^{x_2}(x_1,x_2)-\hat{S}_{x_1}^{x_2}(x_1,x_2))f_1(x_1,x_2).
\end{equation*}
Then, the approximated sensitivity-conditioning system \eqref{eq:sys2predsensapprox} is input-to-state stable with respect to $\xi$.
}
\end{prop}

\begin{proof}
\rr{Since $(x_1,x_2^s(x_1))$ under \eqref{eq:sys2predsens} is exponentially stable under the conditions of Corollary~\ref{cor:sys2stab3}, it is input-to-state stable in \eqref{eq:sys2predsensapprox} with respect to $\xi$ \cite[Lemma~4.6]{khalil2002nonlinear}.
}
\end{proof}
\rr{If the contraction conditions \eqref{eq:sys2condstab} in Corollary~\ref{cor:sys2stab3} do not hold globally, a local result along the line of Proposition~\ref{prop:sys2stab3} can be derived based on the local exponential stability of $x_2^s(x_1)$.}

\section{Example I: Cascade Control}\label{sec:casccont}

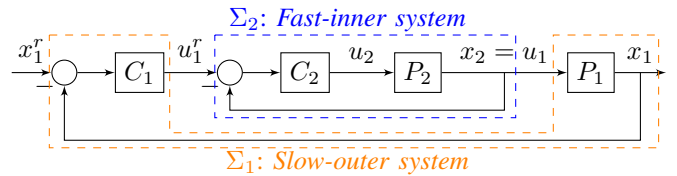
\begin{figure}[b]
    \centering
\begin{tikzpicture}[auto, node distance=2cm,>=latex']
    \node [input, name=inputou] (inputou) {};
    \node [sum, right of=inputou, node distance = 0.7cm] (sum1ou) {};
    \node [block, right of=sum1ou, node distance = 1.cm] (controllerou) {$C_1$};
    
    \node [sum, right of=controllerou, node distance = 1.2cm] (sum1in) {};
    \node [block, right of=sum1in,node distance = 1.cm] (controller) {$C_2$};
    \node [block, right of=controller,node distance=1.5cm] (systemin) {$P_2$};
    \node [block, right of=systemin, node distance=2.3cm] (systemou) {$P_1$};
    \node [tmp, below of=controller,node distance =0.5cm] (tmp1){};
    
    \node [output, right of=systemou, node distance=1.cm] (output) {};
    \node [tmp, below of=controller,node distance = 0.9cm] (tmp2){};
    
    \draw [->] (inputou) -- node{$x_1^r$} (sum1ou);
    \draw [->] (sum1ou) --node[above,name=zou]{} (controllerou);
    \draw [->] (controllerou) -- node[above,name=u1r]{$u_1^r$} (sum1in);
    
    \draw [->] (sum1in) --node[above,name=z]{} (controller);
    \draw [->] (controller) --node[above]{$u_2$} (systemin);
    \draw [->] (systemin) -- node [name=x2]{$x_2=u_1$} (systemou);
    \draw [->] (x2) |- (tmp1)-| node[pos=0.99] {$-$} (sum1in);
    \draw [->] (systemou) --node[above,name=x1]{$x_1$} (output);
    
    \draw [->] (x1) |- (tmp2)-| node[pos=0.99] {$-$} (sum1ou);
    
    \draw [blue,dashed] (2.7,0.5) -- (6.7,0.5) -- (6.7,-0.6) -- (2.7,-0.6) -- (2.7,0.5);
   	\node[text=blue, above of = controller, node distance = 0.7cm] (inner) { \hspace{1cm} $\Sigma_2$: \textit{Fast-inner system}};
   	
   	\draw [orange,dashed] (0.5,0.5) -- (2.1,0.5) -- (2.1,-0.8) -- (7.2,-0.8) -- (7.2,0.5) -- (8.6,0.5) -- (8.6,-1.) -- (0.5,-1.) -- (0.5,0.5) ;
   	\node[text=orange, below of = inner, node distance = 1.9cm] (outer) { \hspace{1cm} $\Sigma_1$: \textit{Slow-outer system}};

\end{tikzpicture}
    \caption{Block diagram of a cascade control.}\label{fig:casccont}
\end{figure}

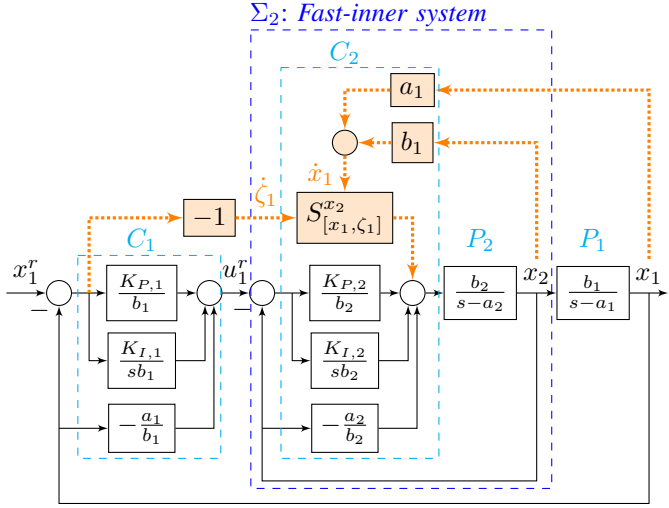
\begin{figure}[t]
\centering
\begin{tikzpicture}[auto, node distance=2cm,>=latex']
    \node [input, name=inputou] (inputou) {};
    \node [sum, right of=inputou, node distance = 0.7cm] (sum1ou) {};
    \node [block, right of=sum1ou, node distance = 1.1cm] (controllerou) {$\frac{K_{P,1}}{b_1}$};
    \node [block, below of=controllerou,node distance=0.9cm] (rateou) {$\frac{K_{I,1}}{sb_1}$};
    \node [block, below of=rateou,node distance=0.9cm] (evolou) {$-\frac{a_1}{b_1}$};
    \node [sum, right of=controllerou,node distance=0.9cm] (sum2ou) {};
    
    \node [sum, right of=sum2ou, node distance = 0.7cm] (sum1in) {};
    \node [block, right of=sum1in,node distance = 1.1cm] (controller) {$\frac{K_{P,2}}{b_2}$};
    \node [block, below of=controller,node distance=0.9cm] (rate) {$\frac{K_{I,2}}{sb_2}$};
    \node [block, below of=rate,node distance=0.9cm] (evolin) {$-\frac{a_2}{b_2}$};
    \node [sum, right of=controller,node distance=0.9cm] (sum2in) {};
    \node [block, right of=sum2in,node distance=0.9cm] (systemin) {$\frac{b_2}{s-a_2}$};
    \node [block, right of=systemin, node distance=1.5cm] (systemou) {$\frac{b_1}{s-a_1}$};
    \node [tmp, below of=controller,node distance =2.5cm] (tmp1){};
    \node [tmp, left of=evolin,node distance =1.1cm] (tmp11){};
    
    \node [output, right of=systemou, node distance=1.cm] (output) {};
    \node [tmp, below of=controller,node distance =2.8cm] (tmp2){};
    \node [tmp, left of=evolou,node distance =1.1cm] (tmp22){};
    
    \draw [->] (inputou) -- node{$x_1^r$} (sum1ou);
    \draw [->] (sum1ou) --node[above,name=zou]{} (controllerou);
    \draw [->] (controllerou) -- (sum2ou);
    \draw [->] (zou) |- (rateou);
    \draw [->] (rateou) -| (sum2ou.250);
    \draw [->] (sum2ou) -- node{$u_1^r$} (sum1in);
    \draw [->] (tmp22) -- (evolou);
    \draw [->] (evolou) -| (sum2ou.290);
    
    \draw [->] (sum1in) --node[above,name=z]{} (controller);
    \draw [->] (controller) -- (sum2in);
    \draw [->] (sum2in) -- node{} (systemin);
    \draw [->] (systemin) -- node [name=x2]{$x_2$} (systemou);
    \draw [->] (z) |- (rate);
    \draw [->] (rate) -| (sum2in.250);
    \draw [->] (x2) |- (tmp1)-| node[pos=0.99] {$-$} (sum1in);
    \draw [->] (systemou) --node[above,name=x1]{$x_1$} (output);
    \draw [->] (tmp11) -- (evolin);
    \draw [->] (evolin) -| (sum2in.290);
    
    \draw [->] (x1) |- (tmp2)-| node[pos=0.99] {$-$} (sum1ou);
    
    \node [block, above of= controller, node distance = 1cm,text = black, fill=orange ,fill opacity=.2, draw opacity=1, text opacity = 1] (sensi) {$S_{[ x_1 , \zeta_1]}^{x_2}$};
    \node[sum, above of = sensi, fill=orange ,fill opacity=.2, draw opacity=1, text opacity = 1] (sumsensi) {}; 
    \node [block, above of= sum2in, node distance = 2cm, text=black, fill=orange ,fill opacity=.2,draw opacity=1, text opacity = 1] (b1) {$b_1$};
    \node [block, above of= sum2in, node distance = 2.7cm, text=black, fill=orange ,fill opacity=.2,draw opacity=1, text opacity = 1] (a1) {$a_1$};
    \node [block, above of= sum2ou, node distance = 1cm, text=black, fill=orange ,fill opacity=.2,draw opacity=1, text opacity = 1] (min1) {$-1$};
    \draw [orange,->,very thick,densely dotted] (x2) |- (b1);
    \draw [orange,->,very thick,densely dotted] (x1) |- (a1);
    \draw [orange,->,very thick,densely dotted] (sumsensi) -- node[left]{$\dot{x}_1$} (sensi);
    \draw [orange,->,very thick,densely dotted] (zou.south) |-  (min1);
    \draw [orange,->,very thick,densely dotted] (min1.east) -- node[above,pos=0.5]{$\dot{\zeta}_1$} (sensi);
    \node [orange,tmp, above of = controller, node distance = 2.7cm] (tmpa1) {}; 
    \draw [orange,->,very thick,densely dotted] (a1.west) |- (tmpa1) -| (sumsensi);
    \draw [orange,->,very thick,densely dotted] (b1.west) -- (sumsensi.east);
    \node [orange,tmp, above of = sum2in, node distance = 1cm] (tmpsum2in) {}; 
    \draw [orange,->,very thick,densely dotted] (sensi.east) |- (tmpsum2in)-| (sum2in);
    
    \draw [blue,dashed] (3.25,3.5) -- (7.25,3.5) -- (7.25,-2.6) -- (3.25,-2.6) -- (3.25,3.5);
     \node[text=blue, above of = controller, node distance =3.7cm] (inner) {\hspace{0.5cm} $\Sigma_2$: \textit{Fast-inner system}};
     \draw [cyan,dashed] (3.65,3) -- (5.75,3) -- (5.75,-2.2) -- (3.65,-2.2) -- (3.65,3);
     \draw [cyan,dashed] (0.95,0.5) -- (2.85,0.5) -- (2.85,-2.2) -- (0.95,-2.2) -- (0.95,0.5);
     \node[text=cyan, above of = controller, node distance = 3.2cm] {$C_2$};
     \node[text=cyan, above of = controllerou, node distance = 0.7cm] {$C_1$};
     \node[text=cyan, above of = systemin, node distance = 0.7cm] {$P_2$};
     \node[text=cyan, above of = systemou, node distance = 0.7cm] {$P_1$};          
    
\end{tikzpicture}
\caption{Block diagram of cascade PI control \eqref{eq:cascPI}, with plants $P_i$ and controllers $C_i$ in Laplace domain. The additional sensitivity-conditioning elements in \eqref{eq:cascPIpred} are depicted in orange blocks and dotted thicker arrows.}\label{fig:cascPI}
\end{figure}

Consider a standard cascade control architecture \cite{lee1998pidcascade}, see Fig.~\ref{fig:casccont}, with a fast-inner closep-loop system $\Sigma_2$, and a slow-outer closed-loop system $\Sigma_1$, both with plants $P_i$ and controllers $C_i$. 
More concretely, consider an example with two linear scalar first-order systems: 
\begin{equation*}\label{eq:cascplants}\begin{array}{rl}
    P_1: & \dot{x}_1 = a_1 x_1 + b_1 u_1, \; u_1 =  x_2  \\
    P_2: & \dot{x}_2 = a_{2} x_2 + b_2 u_2,  
\end{array}
\end{equation*}
where all parameters are real-valued. PI (proportional-integral) controllers are typically used for $C_1,C_2$: 
\begin{equation}\label{eq:cascPIcont}\begin{array}{rl}
     C_1: & x_2^r = u_1^r = -\frac{1}{b_1}\big( a_1 x_1 +K_{P,1}(x_1-x_1^r) +K_{I,1}\zeta_1 \big) \\  
     & \dot{\zeta}_1 = (x_1-x_1^r) \\
     C_2: & u_2 = -\frac{1}{b_2}\big(a_{2} x_2 +K_{P,2}(x_2-x_2^r) +K_{I,2}\zeta_2 \big) \\ 
     & \dot{\zeta}_2 = (x_2-x_2^r).
\end{array}
\end{equation}
where $\zeta_i$ are the integral error states, $K_{P,i},K_{I,i}$ are control gains to be determined, and the terms $a_i x_i$ are feed-forward terms to cancel the system dynamics. If the systems $\Sigma_i$ had a time-scale separation as \eqref{eq:sys2ideal}, the resulting interconnected system, with states $x_i,\zeta_i$ for each $\Sigma_i$, can be expressed as:
\begin{equation}\label{eq:cascPI}\begin{array}{rcl}
    \Sigma_1: & \dot{x}_1 = & -K_{P,1}(x_1-x_1^r) -K_{I,1}\zeta_1\\
    & \dot{\zeta}_1 = & (x_1-x_1^r) \\[0.1cm]
    \Sigma_2: & \frac{d x_2}{d\tau} = & -K_{P,2}(x_2-x_2^r) -K_{I,2}\zeta_2\\
    & \frac{d \zeta_2}{d\tau} = & (x_2-x_2^r) \\
    & x_2^r = & -\frac{1}{b_1}\big( a_1 x_1 +K_{P,1}(x_1-x_1^r) +K_{I,1}\zeta_1 \big)
\end{array}
\end{equation}
which admits the globally asymptotically stable steady state $x_1^s=x_1^r,x_2^s=x_2^r=u_1^r=\frac{-a_1 x_1^r}{b_1},\zeta_1^s=0,\zeta^s=0$ for positive gains $K_{P,i}>0, K_{I,i}>0$. See Fig.~\ref{fig:cascPI} for a block diagram representation of this control architecture. 

\begin{rem}
The feed-forward control inputs in \eqref{eq:cascPIcont} can also be implemented using the references $x_1^r,x_2^r(=u_1^r)$ instead of the states $x_1,x_2$, to compensate the plant dynamics. Then, the conditions for asymptotic stability of \eqref{eq:cascPI} are $K_{P,i}>a_i$ and $K_{I,i}>0$. These controllers \eqref{eq:cascPIcont} may also not include any feed-forward compensation at all. Then, the conditions for asymptotic stability of \eqref{eq:cascPI} are $a_i-b_i K_{P,i}<0, K_{I,i}>0$. In either case, all our subsequent results hold with minor adjustments.
\end{rem}

\rr{To preserve the stability of the time-scale separated cascaded system with controllers \eqref{eq:cascPI} in a single time scale we apply the sensitivity-conditioning \eqref{eq:sys2predsens}:}
\begin{equation}\label{eq:cascPIpred}\begin{array}{l}
     \dot{x}_1 =   a_1x_1 +b_1 x_2  \\
     \dot{\zeta}_1 =  (x_1-x_1^r) \\
     x_2^r = u_1^r =  -\frac{1}{b_1}\big( a_1 x_1 +K_{P,1}(x_1-x_1^r) +K_{I,1}\zeta_1 \big) \\
     \dot{x}_2 =  -K_{P,2}(x_2-x_2^r) -K_{I,2}\zeta_2  + S_{[x_1 ,\zeta_1]}^{x_2} 
     \begin{bmatrix} \dot{x}_1 \\ \dot{\zeta}_1 \end{bmatrix}\\[-0.3cm]
     \dot{\zeta}_2 =  (x_2-x_2^r) + S_{[x_1 ,\zeta_1]}^{\zeta_2} \begin{bmatrix} \dot{x}_1 \\ \dot{\zeta}_1 \end{bmatrix} \hspace{-0.1cm}, \\
     \end{array}
\end{equation}
where
\begin{equation*}
     \begin{array}{rl}
        \begin{bmatrix} S_{[x_1 ,\zeta_1]}^{x_2} \\ S_{[x_1 ,\zeta_1]}^{\zeta_2} \end{bmatrix} = & 
        - \hspace{-0.1cm} \begin{bmatrix} -K_{P,2} & -K_{I,2} \\ 1 & 0 \end{bmatrix}^{-1} \hspace{-0.15cm}
     \begin{bmatrix} K_{P,2} \\ -1 \end{bmatrix} \hspace{-0.2cm}
     \begin{bmatrix}-\frac{a_1+K_{P,1}}{b_1} \\ -\frac{K_{I,1}}{b_1} \end{bmatrix}^T \\
         = & \begin{bmatrix}-\frac{a_1+K_{P,1}}{b_1} & -\frac{K_{I,1}}{b_1} \\ 0 & 0\end{bmatrix} \hspace{-0.1cm}.
     \end{array} 
\end{equation*}
Note that $S_{[x_1 , \zeta_1]}^{\zeta_2}=[0 \; 0]$ is due to $\zeta_2^s=0$ for all $x_1,\zeta_1$, for a stable inner system $\Sigma_2$. In summary, the control structure can be graphically represented as in Fig.~\ref{fig:cascPI}, with the sensitivity-conditioning elements acting as a derivative-type control.

In compact matrix form, the closed-loop system reads as
\begin{equation*}\begin{array}{l}
    \dot{x} = T A x + B x_1^r, \; x = [ {x_1}^T,{\zeta_1}^T,{x_2}^T,{\zeta_2}^T ]^T \\[0.1cm]
     A = \begin{bmatrix} a_1 & 0 & b_1 & 0 \\ 1 & 0 & 0 & 0 \\
    -K_{P,2}\frac{a_1+K_{P,1}}{b_1} & -K_{P,2}\frac{K_{I,1}}{b_1} & -K_{P,2} & -K_{I,2} \\
    \frac{a_1+K_{P,1}}{b_1} & \frac{K_{I,1}}{b_1} & 1 & 0
    \end{bmatrix}  \\[0.9cm]
     T = \begin{bmatrix} I & 0 \\ S_{[x_1 ,\zeta_1]}^{[x_2,\zeta_2]} & I \end{bmatrix}, 
     \; T^{-1} =  \begin{bmatrix} I & 0 \\ -S_{[x_1 ,\zeta_1]}^{[x_2,\zeta_2]} & I \end{bmatrix} \\
     B = \begin{bmatrix} 0 \; & \; -1  \; & \; K_{P,2}\frac{K_{P,1}}{b_1} \; & \; -\frac{K_{P,1}}{b_1} \end{bmatrix}^T \hspace{-0.1cm}.
\end{array}
\end{equation*}
By means of the similarity transformation used in the proof of Proposition~\ref{prop:sys2stab2}, we obtain
\begin{equation*}
    TA \sim T^{-1}(TA)T=\begin{bmatrix} -K_{P,1} & -K_{I,1} & \star & \star \\ 
    1 & 0 & \star & \star \\
    0 & 0 & -K_{P,2} & -K_{I,2} \\
    0 & 0 & 1 & 0
    \end{bmatrix} \hspace{-0.1cm},
\end{equation*}
where $\star$ are irrelevant terms for the following considerations. This block-companion form of $TA$ confirms that the only conditions required for stability of \eqref{eq:cascPIpred} are $K_{P,i}>0$ and $K_{I,i}>0$, as for \eqref{eq:cascPI}. This is aligned with Proposition~\ref{prop:sys2stab2}: stability of \eqref{eq:cascPI} is preserved in \eqref{eq:cascPIpred} using the sensitivity-conditioning \eqref{eq:sys2predsens}. Note that the system \eqref{eq:cascPIpred} without the sensitivity-conditioning term would have system matrix $A$ instead of $TA$: $\dot{x} = A x + B x_1^r$. Then, for example with $a_i=0,b_i=1,K_{P,i}=1>0,K_{I,i}=1>0$, $A$ has positive eigenvalues despite having positive control parameters, so stability is lost without the sensitivity-conditioning term.

\begin{rem}\label{rem:backstepcasc}
\rr{In cascade control the steady-state closed form $x_2^s(x_1)$ in \eqref{eq:sys2ideal} is typically known by design. Hence, we could use an alternative sensitivity based on the Lyapunov function $\normsz{x_2-x_2^s(x_1)}^2_{P_2}$, see Remark~\ref{rem:backstep}. Therefore, in this cascade control example \eqref{eq:cascPI}, the sensitivity-conditioning approach \eqref{eq:sys2predsens} is equivalent to backstepping \cite[Ch.~14]{khalil2002nonlinear} for nonlinear control design, up to an extra proportional control term.} 
\end{rem}

\subsection{Numerical simulation: DC/AC-converter + RLC filter}

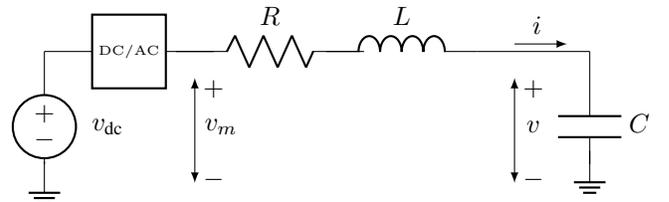
\begin{figure}[b]
    \centering
 \begin{tikzpicture}
\draw (-0.750,0) to [twoport, t=$\scriptscriptstyle \mathrm{DC/AC}$](1.5,0);
\draw (1,0) [R=$R$] to (3.5,0);
\draw (3,0)  [L=$L$] to (5,0);
\draw (5,0) -- (6.5,0);
\draw (6.5,0) to (6.5,-0.5) [C=$C$]  to node[ground]{} (6.5,-1.5);

\draw (5.5,0.1) [- latex] to (6.2,0.1);
\node at(5.8,0.35) {$i$};

\draw[latex - latex]  (5.5,-0.35) to (5.5,-1.75);
\node at (5.75,-0.5) {$+$};
\node at (5.75,-1.7) {$-$};
\node at (5.75,-1.) {$v$};

\draw (-0.75,0) to (-0.75,-0.50)[american voltage source=$v_{\text{\upshape{dc}}}$]  to node[ground]{} (-0.75,-1.5);

\draw [latex - latex] (1.25,-0.35)  to (1.25,-1.75);
\node at (1.5,-0.5) {$+$};
\node at (1.5,-1.7) {$-$};
\node  at (1.6,-1.) {$v_{m}$};
\end{tikzpicture}
    \caption{DC/AC-converter with RLC circuit}
    \label{fig:rlccirc}
\end{figure}

Designing a DC/AC-converter connected to a $RLC$-filter \cite{yazdani2010voltage} is a standard control problem in power electronics. Here we present a simplified version, where the DC/AC-converter modulates the DC voltage $v_{\text{\upshape{dc}}}$ into the three-phase AC voltage $v_m$. Using an averaged converter with stiff DC voltage, $v_m$ is a fully controllable voltage source. This modulated voltage $v_m$ is then used to control the three-phase current $i$ through the resistance $R$ and inductance $L$, which in turn is used to control the output voltage $v$ at the capacitor $C$ to follow a reference $ v^r$, see Fig.~\ref{fig:rlccirc}. 


Let the electrical signals be represented in rectangular coordinates, using the real and imaginary parts, so $i,v,v_m \in \mathbb{R}^2$, and define the frequency $\omega$ and the rotation matrix $\mathcal{J}=\left[\begin{smallmatrix}0 &-1 \\ 1 & 0 \end{smallmatrix}\right]$. Then, according to Kirchoff's laws, the electrical signals dynamics in the rotation frame coordinates \cite{subotic2019lyapunov} are:
\begin{equation*}
    C\frac{dv}{dt} = i - \mathcal{J}\omega C v, \; L\frac{di}{dt} = v_m -(R+\mathcal{J} \omega L)i - v,
\end{equation*}

Following the cascaded PI example \eqref{eq:cascPI}, the reference $i^r$ and the controller $v_m$ are chosen as PI controllers with every $K_{(\cdot,\cdot)}=k_{(\cdot,\cdot)}I_2$, where $I_2$ is identity of dimension $2$:
\begin{equation}\label{eq:cascPIex}\begin{array}{rl}
     \frac{dv}{dt} = & C^{-1} i - \mathcal{J} \omega v \\
     \frac{d\zeta_v}{dt} = & (v-v^r) \\
     i^r = & \mathcal{J}\omega C v + C(-K_{P,v}(v-v^r) -K_{I,v}\zeta_v ) \\[0.1cm]
     v_m = & (R+\mathcal{J} \omega L)i+v +L(-K_{P,i}(i-i^r) -K_{I,i}\zeta_i) \\
     \frac{di}{dt} = & -K_{P,i}(i-i^r) -K_{I,i}\zeta_i \\
     \frac{d\zeta_i}{dt} = & (i-i^r)
\end{array}
\end{equation}

\begin{table}[t]
    \centering
    \caption{Simulation parameters}
    \label{tab:param}
\begin{tabular}{|l|l|}
\hline
    RLC-filter & $R=1m\Omega, L= 1mH, C=300\mu F$ \\[0.05cm] \hline
    Frequency & $f = 50 \frac{1}{\text{sec}}, \; \omega=2 \pi 50 \frac{\text{rad}}{\text{sec}}$ \\[0.05cm] \hline
    Outer Controller $C_1$ & $k_{P,v}=30 \frac{A}{VF}, k_{I,v}=0.3\frac{A}{VF}\frac{\text{rad}}{\text{sec}}$ \\[0.05cm] \hline
    Reference & Real and imaginary parts: $v_\Re^r=120V, \;v_\Im^r=0V$ \\
    & Magnitude: $|v^r|=120V=1$ p.u. (per unit)  \\[0.05cm] \hline
    Black start & $v(0)=0V=0$ p.u., $i(0)=0A=0$ p.u. \\[0.05cm] \hline
\end{tabular}
\end{table}


\begin{figure}[t]
    \centering
    \includegraphics[width=9.cm]{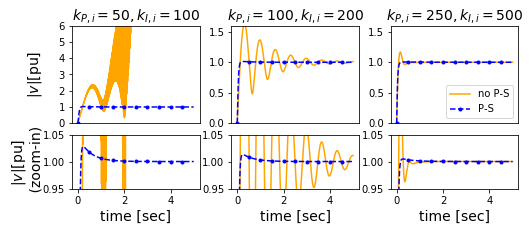}
    \caption{Simulation of a black-start of the electrical system \eqref{eq:cascPIex} with increasing values for the control parameters $k_{P,i},k_{I,i}$, with and without sensitivity-conditioning (P-S).}
    \label{fig:cascPIex}
\end{figure}

Fig.~\ref{fig:cascPIex} shows simulation results for \eqref{eq:cascPIex} with the parameters in Table~\ref{tab:param}. We simulate a black start, i.e., the system starts at time $t=0$ with zero values, and needs to track a given constant reference, which can be interpreted as a step input. \rr{Since system \eqref{eq:cascPIex} does not have an intrinsic time-scale separation, we use sufficiently large gains $k_{P,i},k_{I,i}$ for the inner controller $C_2$ to artificially enforce a wider time-scale separation}, i.e., lower $\epsilon$ in \eqref{eq:sys2singpert}:
\begin{itemize}[leftmargin=*]
    \item $k_{P,i}=50 \frac{V}{AH}, k_{I,i}=100\frac{V}{AH}\frac{\text{rad}}{\text{sec}}$: the system without sensitivity-conditioning is unstable even though the control parameters are all positive. On the contrary, the sensitivity-conditioning turns it into a stable system, with the magnitude $|v|$ stabilizing quickly at $1$ p.u. and the frequency $f$ at $50\frac{1}{\text{sec}}$.
    \item $k_{P,i}=100\frac{V}{AH}, k_{I,i}=200\frac{V}{AH}\frac{\text{rad}}{\text{sec}}$: the system without sensitivity-conditioning becomes stable, but it is still an unacceptable controller due to the high overshoot and relatively large settling time. Again the sensitivity-conditioning turns it into a well-performing controller.
    \item $k_{P,i}=250\frac{V}{AH}, k_{I,i}=500\frac{V}{AH}\frac{\text{rad}}{\text{sec}}$: the system without sensitivity-conditioning becomes stable with acceptable control performance. However, the sensitivity-conditioning approach performs much better with negligible overshoots.
\end{itemize}
Moreover, note in the zoom-in in Fig.~\ref{fig:cascPIex} that the overshoot decreases as $k_{P,i},k_{I,i}$ increase, so the system with sensitivity-conditioning also benefits from having a faster controller $C_2$. Such faster $C_2$ can further increase the convergence rate.

\section{Example II: Bilevel Optimization}\label{sec:bilopt}

\rr{In this section, we show an application of the sensitivity-conditioning \eqref{eq:sys2predsens} to bilevel optimization \cite{bard1998practical}. As opposed to the cascade control example in Section~\ref{sec:casccont}, in this case the steady-state map $x_2^s(x_1)$ is not available in closed form, thus backstepping is not applicable, see Remarks~\ref{rem:backstep} and \ref{rem:backstepcasc}. Yet, the sensitivity-conditioning \eqref{eq:sys2predsens} can still be used to preserve the stability of the two-time-scale system \eqref{eq:sys2ideal} in a single one.}

Consider a general unconstrained bilevel problem \cite{sinha2018reviewbil, bard1998practical}:

\begin{equation}\label{eq:bilopt}\begin{array}{ll}
    \min_{x_1,x_2^*} F_1(x_1,x_2^*) \\[0.0cm]
    \text{ s.t. } x_2^* \in \arg\min_{x_2} F_2(x_1,x_2)
\end{array}
\end{equation}
where $F_1(\cdot),F_2(\cdot)$ are the upper- and lower-level objective functions, respectively. 

\begin{ass}[Adaptation of Assumptions~\ref{ass:sys2sensinv} and \ref{ass:sys2loclip}]\label{ass:bilsensinv}
\
\begin{itemize}[leftmargin=*]
    \item The functions $F_1(\cdot)$ and $F_2(\cdot)$ are twice and thrice continuously differentiable, respectively, with Lipschitz continuous partial derivatives.
    \item \rr{For every $x_1$, the lower-level problem $\arg\min_{x_2} F_2(x_1,x_2)$ has at most a single solution $x_2^*$, where the second-order partial derivative $\nabla_{x_2x_2}^2 f_2(x_1,x_2^*)$ is invertible.}
\end{itemize}
\end{ass}
\rr{The single solution assumption is often used as simplification in bilevel problems \cite{bard1998practical}. It ensures that \eqref{eq:bilopt} is well-posed, and it allows to simplify the constraint to $x_2^* = \arg\min_{x_2} F_2(x_1,x_2)$.} Moreover, the invertibility assumption allows to define the sensitivity of $x_2^*$, a known concept in bilevel optimization \cite{dempe2012sensitivity}, similar to the one introduced in \eqref{eq:sys2sens}: Consider a point $(x_1,x_2^*)$ satisfying the first-order optimality conditions of the lower-level problem, i.e., $\nabla_{x_2} F_2(x_1,x_2^*)=0$. Since $\nabla_{x_2 x_2}^2 F_2(x_1,x_2^*)$ is invertible, the implicit function theorem \cite{krantz2012implicit} guarantees the local existence of the map $x_2^*(x_1)$, and gives an expression for its derivative:
\begin{equation*}\begin{array}{l}
     \nabla_{x_1} x_2^*(x_1)  \hspace{-0.05cm} = \hspace{-0.05cm} -\big( \hspace{-0.02cm} \nabla_{x_2 x_2}^2 \hspace{-0.05cm} F_2(x_1,x_2^*(x_1) \hspace{-0.03cm} ) \hspace{-0.02cm} \big)^{\hspace{-0.05cm} -1}\nabla_{x_2 x_1}^2\hspace{-0.05cm}  F_2(x_1,x_2^*(x_1)\hspace{-0.03cm}) 
\end{array}
\end{equation*}

Additionally, $x_2^*(x_1)$ can be used to locally define a reduced objective $F_1$: $F_1^r(x_1):=F_1(x_1,x_2^*(x_1))$, as in \eqref{eq:sys2ideal}, and use $\nabla_{x_1} x_2^*(x_1)$ to give an expression for the \textit{total derivative},
\begin{equation}\label{eq:biltotder} \begin{array}{l}
     D_{x_1} F_1(x_1,x_2^*(x_1)) := \nabla_{x_1} F_1^r(x_1) \\
     = \nabla_{x_1} F_1(x_1,x_2^*(x_1)) + \nabla_{x_1} x_2^*(x_1)^T \nabla_{x_2}F_2(x_1,x_2^*(x_1))
\end{array} 
\end{equation}
defined for points where $x_2=x_2^*(x_1)$. Under Assumption~\ref{ass:bilsensinv}, $\nabla_{x_2x_2}^2 f_2(x_1,x_2)$ is invertible for $x_2$ in a neighborhood of $x_2^s(x_1)$. Hence, $\nabla_{x_1} x_2^*(x_1)$ and $\nabla_{x_1} F_1^r(x_1)$ can be extended to these $x_2$, as the extended sensitivity in \eqref{eq:sys2sens}:
\begin{equation}\label{eq:bilsenstotderext} \begin{array}{rl}
 D_{x_1} \hspace{-0.05cm} F_1(x_1,x_2) \hspace{-0.05cm} & \hspace{-0.05cm} := \hspace{-0.05cm} \nabla_{x_1} \hspace{-0.05cm} F_1(x_1,x_2) \hspace{-0.07cm} + \hspace{-0.07cm} S_{x_1}^{x_2}(x_1,x_2)^T \nabla_{x_2} \hspace{-0.05cm} F_1(x_1,x_2) \\
S_{x_1}^{x_2}(x_1,x_2) \hspace{-0.05cm} & \hspace{-0.05cm} := \hspace{-0.05cm} -(\nabla_{x_2 x_2}^2 F_2(x_1,x_2))^{-1}\nabla_{x_2 x_1}^2 F_2(x_1,x_2),
 \end{array}
\end{equation}
satisfying the restrictions $D_{x_1} F_1(x_1,x_2)_{|_{(x_1,x_2^*(x_1))}} \allowbreak =\nabla_{x_1} F_1^r(x_1)$, $S_{x_1}^{x_2} (x_1,x_2)_{|_{(x_1,x_2^*(x_1))}}=\nabla_{x_1} x_2^*(x_1)$.



\subsection{Bilevel local solutions}
Understanding the properties of the bilevel problem solutions is essential to connect the convergence of algorithms with the stability of steady states from previous sections.
Therefore, we recall the concept of \textit{local solutions} in \cite[Ch.~8]{bard1998practical} to represent locals minima of \eqref{eq:bilopt} and their first and second-order necessary and sufficient conditions:

\begin{defi}[local solution] \label{def:localbil} \cite[Ch.~8]{bard1998practical}
A point $(x_1^*,x_2^*)$ is a \textit{(strict) local solution} of \eqref{eq:bilopt} if:
\begin{enumerate}[leftmargin=*]
    \item The point $x_2^*$ is a local minimum of \rr{$F_2(x_1^*,\cdot)$ with fixed $x_1^*$}. 
    \item There exists a neighborhood $N$ of $(x_1^*,x_2^*)$ such that $F_1(x_1^*,x_2^*) \leq F_1(x_1,x_2)$ ($<$ for strict) for all $(x_1,x_2)\in N$ such that $x_2$ is a local minimum of \rr{$F_2(x_1,\cdot)$ with fixed $x_1$}. 
\end{enumerate}
\end{defi}

\begin{prop}[First-order necessary conditions]\label{prop:critpoint} \cite[Ch.~8]{bard1998practical}
A \textit{local solution} $(x_1^*,x_2^*)$ is a \textit{stationary point}, i.e., it satisfies 
$\nabla_{x_2} F_2(x_1^*,x_2^*)=0, \; \nabla_{x_1} F_1^r(x_1^*) =0.$
\end{prop}
\begin{prop}[Second-order conditions] \label{prop:cond2ord} \cite[Ch.~8]{bard1998practical}
\begin{itemize}[leftmargin=*]
\item \textit{Necessary conditions}: A \textit{local solution} $(x_1^*,x_2^*)$ satisfies
\begin{equation*}\label{eq:cond2ord}\begin{array}{l}
    \nabla_{x_2 x_2}^2 F_2(x_1^*,x_2^*) \succeq 0 , \; 
    \nabla_{x_1 x_1}^2 F_1^r(x_1^*) \succeq 0 
\end{array}
\end{equation*}
\item \textit{Sufficient conditions}: A \textit{stationary point} $(x_1^*,x_2^*)$ satisfying
\begin{equation}\label{eq:cond2ordstr}\begin{array}{l}
    \nabla_{x_2 x_2}^2 F_2(x_1^*,x_2^*) \succ 0 , \; 
     \nabla_{x_1 x_1}^2 F_1^r(x_1^*) \succ 0 ,
\end{array}
\end{equation}
is a \textit{strict local solution}.
\end{itemize}
\end{prop}

%

\subsection{Bilevel gradient flow}
The steepest descent direction method \cite{savard1994steepbil} is a standard approach to iteratively solve \eqref{eq:bilopt}. It follows the negative gradient of $F_1^r(x_1)$ 
with step size $\alpha^k$ in each iteration $k$:
\begin{equation}\label{eq:bilgradu}\begin{array}{rl}
    x_1^{k+1} = & x_1^k - \alpha^k \nabla_{x_1} F_1^r(x_1^k) = x_1^k - \alpha^k \big(D_{x_1} F_1(x_1^k,x_2^k) \big) \\
        x_2^k = & \arg\min_{x_2} F_2(x_1^{k},x_2)
\end{array}
\end{equation}
If the lower-level update $x_2^k = \arg\min_{x_2} F_2(x_1^{k},x_2)$ is not available in closed form, it can be solved iteratively using for example gradient descent with step size $\beta^l$ and updates: $${x}_2^{l+1} = {x}_2^{l} - \beta^l \nabla_{x_2} F_2(x_1^k,{x}_2^{l})$$ 
The corresponding continuous-time version of this bilevel gradient descent \eqref{eq:bilgradu} can be represented on two time scales with $\epsilon \to 0$ and the singular perturbation interconnection \eqref{eq:sys2singpert}:
 
\begin{equation}\label{eq:bilgradfasyn}\begin{array}{rl}
    \dot{x}_1 = & -D_{x_1} F_1(x_1,x_2) \\
    \epsilon\dot{x}_2 = & -\nabla_{x_2} F_2(x_1,x_2)
\end{array}
\end{equation} 

As mentioned before in Section~\ref{sec:moti}, these nested iterations \eqref{eq:bilgradu} on two time scales \eqref{eq:bilgradfasyn} may slow down the algorithm convergence. On the other hand, 
the sensitivity-conditioning system \eqref{eq:sys2predsens} yields:

\begin{equation}\label{eq:bilgradfcorr}\begin{array}{rl}
    \dot{x}_1 = & -D_{x_1} F_1(x_1,x_2) \\
    \dot{x}_2 = & -\nabla_{x_2} F_2(x_1,x_2) + S_{x_1}^{x_2}(x_1,x_2) \; \dot{x}_1
\end{array}
\end{equation}

\begin{cor}[Local convergence of \eqref{eq:bilgradfcorr}]\label{cor:bilstab}
A point $(x_1^*,x_2^*)$ is a \textit{strict local solution} of \eqref{eq:bilopt} satisfying the sufficient conditions in \eqref{eq:cond2ordstr} if and only if it is a locally exponentially stable steady state of the sensitivity-conditioning bilevel gradient flow \eqref{eq:bilgradfcorr}.
\end{cor}

\begin{proof}
Note that Assumption~\ref{ass:bilsensinv} adapts Assumptions~\ref{ass:sys2sensinv} and \ref{ass:sys2loclip} for the bilevel problem \eqref{eq:bilopt}, and the second-order total derivative 
is symmetric and satisfies  
$$ \begin{array}{l}
D_{x_1 x_1}^2 F_1(x_1,x_2^*(x_1)) := \nabla_{x_1 x_1}^2 F_1^r(x_1) = \\
\nabla_{x_1} \hspace{-0.03cm} D_{x_1} \hspace{-0.03cm} F_1(x_1,x_2^*(x_1) \hspace{-0.03cm})  
\hspace{-0.07cm} + \hspace{-0.07cm} \nabla_{x_1}{x_2}^*\hspace{-0.03cm}(x_1\hspace{-0.03cm})^T \nabla_{x_2} \hspace{-0.03cm} D_{x_1} F_1(x_1,x_2^*(x_1)\hspace{-0.03cm} )
\end{array}$$ 
Hence, $(x_1^*,x_2^*)$ is locally exponentially stable if and only if \eqref{eq:cond2ordstr} holds \cite[Cor.~4.3]{khalil2002nonlinear}.
\end{proof}

\begin{rem}\label{rem:bilmultisol}
This convergence result can be stated to larger regions if similar conditions as in Proposition~\ref{prop:sys2stab3} hold. 
\end{rem}




\subsection{Time discretization and numerical simulation}

The Euler-forward method \cite{atkinson2008introduction} with time constant $\tau$ can be used to integrate the differential equations \eqref{eq:bilgradfcorr} and \eqref{eq:bilgradfasyn} for a fixed $\epsilon$. Then we get a discrete-time descent algorithm:

\begin{equation}\label{eq:bildisc}\begin{array}{rl}
    x_1^{k+1} = & x_1^k -\tau D_{x_1} F_1(x_1^k,x_2^k) \\
    \text{} \eqref{eq:bilgradfasyn}: \; x_2^{k+1} = & x_2^k -\frac{\tau}{\epsilon }\nabla_{x_2} F_2(x_1^k,x_2^k) \\
    \text{} \eqref{eq:bilgradfcorr}: \; x_2^{k+1} = & x_2^k - \tau \nabla_{x_2} F_2(x_1^k,x_2^k) \\
    & + S_{x_1}^{x_2}(x_1^k,x_2^k) (-\tau D_{x_1} F_1(x_1^k,x_2^k)),
\end{array}
\end{equation}
where the time constant $\tau$ plays the role of the step size in optimization \cite{bertsekas1997nonlinear}. If \eqref{eq:bilgradfcorr} is locally or globally exponentially stable, see Proposition~\ref{prop:sys2stab2} and \ref{prop:sys2stab3}, then under some conditions its Euler-forward discretization will retain this exponential stability for suitable time constants below a certain threshold $\tau < \bar{\tau}$ \cite{stetter1973analysis,stuart1994numerical}. More concretely, this can be proven by extending \cite[Lemma~5]{qu2019expstabsaddle} to the Krasovskii Lyapunov functions in Proposition~\ref{prop:sys2stab3}. 
We will formalize further results for discrete-time sensitivity-conditioning in Subsection~\ref{sec:predsensdisc}.

\begin{rem}
Bilevel optimization problems like \eqref{eq:bilopt} can also be represented as Stackelberg games \cite{von2010market}. The singular perturbation dynamics \eqref{eq:bilgradfasyn}, with its discrete-time version in \eqref{eq:bildisc}, can be interpreted as a simultaneous gradient descent on both variables with different step sizes: $\tau, \frac{\tau}{\epsilon }$. This corresponds to deterministic Stackelberg learning dynamics \cite{fiez2019convergence, fiez2020implicit}. 
\newline
The particular case when $F_2(x_1,x_2)=-F_1(x_1,x_2)$ in \eqref{eq:bilopt}, is called a zero-sum or minimax game \cite{jin2019local, fiez2019convergence}. In this context, this simultaneous gradient descent algorithm is known as the $\gamma$-gradient descent ascent ($\gamma$-GDA) \cite{jin2019local}, with $\gamma=\frac{1}{\epsilon}$. Then, the discrete-time sensitivity-conditioning application for bilevel optimziation in \eqref{eq:bildisc}, corresponds to the Stackelberg generalization of the algorithm in \cite{wang2019followridge} for minimax games. 
\end{rem}

\begin{figure}[b]
    \centering
    \includegraphics[width=9cm]{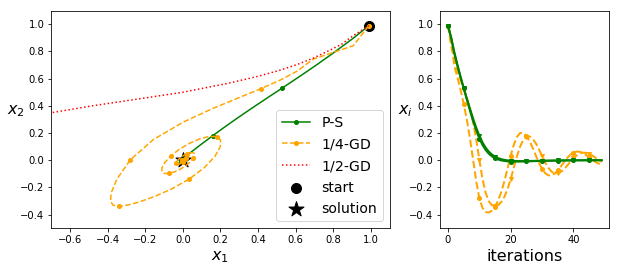}
    \caption{Comparison of discrete-time methods \eqref{eq:bildisc} for the functions in \eqref{eq:exbilscaled}: a simultaneous gradient descent \eqref{eq:bilgradfasyn} for two values $\epsilon \in \{\frac{1}{4},\frac{1}{2}\}$ ($\epsilon$-GD), against the sensitivity-conditioning (P-S) approach \eqref{eq:bilgradfcorr}. Both are implemented in discrete time with markers every $5$ iterations.}
    \label{fig:bilopt}
\end{figure}

Consider an example similar to the ones in \cite{fiez2019convergence, wang2019followridge}: 
\begin{equation}\label{eq:exbilscaled}
F_1(x_1,x_2) = -\frac{x_1^2}{2}+x_2^2, \; F_2(x_1,x_2) = \Big( \frac{x_2^2}{4}-\frac{x_1 x_2}{2} \Big) e^{- \frac{x_2^2}{2}}
\end{equation}
The point $(0,0)$ is a \textit{strict local solution} of the bilevel problem \eqref{eq:bilopt}, since $\nabla_{x_2} \hspace{-0.05cm} F_2(0,0) \hspace{-0.05cm} = \hspace{-0.05cm} 0, D_{x_1} \hspace{-0.05cm} F_1(0,0) \hspace{-0.05cm} = \hspace{-0.05cm} 0, \allowbreak \nabla_{x_2 x_2}^2 F_2(0,0)\hspace{-0.05cm}  =\hspace{-0.05cm} \frac{1}{2} \hspace{-0.05cm} > \hspace{-0.05cm} 0, D_{x_1 x_1}^2 F_1(0,0) \hspace{-0.05cm} = \hspace{-0.05cm} 1>0$. 

In Fig.~\ref{fig:bilopt}, we compare the approaches in \eqref{eq:bildisc} using $\tau=\frac{1}{4}$, singular perturbation values $\epsilon \in \{\frac{1}{4},\frac{1}{2}\}$, and the functions in \eqref{eq:exbilscaled}. First, we observe that the sensitivity-conditioning (P-S) approach \eqref{eq:bilgradfcorr} is able to converge quickly to the solution $(0,0)$. For the case $\epsilon=\frac{1}{2}$, the simultaneous gradient descent ($\epsilon$-GD) based on \eqref{eq:bilgradfasyn} fails to converge. It converges for $\epsilon=\frac{1}{4}$, but still oscillates around the solution causing a slower convergence. 

\section{Sensitivity-conditioning for Multiple Systems}\label{sec:multits}

The theoretical results in Section~\ref{sec:predsens2sys}, and the applications to cascade control and bilevel optimization, Sections~\ref{sec:casccont} and \ref{sec:bilopt} respectively, deal with two-time-scale systems. In this section we show how the sensitivity-conditioning interconnection \eqref{eq:sys2predsens} can be extended to multiple time scales arising in, e.g., multiple nested systems \cite{grammel2004nonlinear} or multilevel programming \cite{BARD1982multilevel}. Consider $N$ {differential-algebraic-equation} subsystems $\Sigma_i$, with states $x_i \in \mathbb{R}^{n_i}$ and vector fields $f_i(\cdot)$, operating on different time scales $\tau_i$ as in \eqref{eq:sys2ideal0}, ordered from slow $\Sigma_1$ to fast $\Sigma_N$:
\begin{equation}\label{eq:sysN}\begin{array}{rl}
    \Sigma_i: \frac{dx_i}{d\tau_i} = f_i(x_1,\dots,x_N) \text{ s.t. } & f_j(x_1,\dots,x_N) =0 \; \forall j>i\\ 
     & \frac{dx_j}{d\tau_i} =0 \; \forall j<i
\end{array}
\end{equation}
As in \eqref{eq:sys2singpert}, consider the corresponding singular-perturbed system with terms $0< \epsilon_N \ll \dots \ll \epsilon_2 \ll \epsilon_1 = 1$ \cite{grammel2004nonlinear}:

\begin{equation*}\label{eq:sysNsingpert}
    \begin{bmatrix} \epsilon_1 I & 0 & 0 \\[-0.15cm]
    0 & \ddots & 0 \\
    0 & 0 & \epsilon_N I \end{bmatrix} \hspace{-0.2cm}
    \begin{bmatrix} \dot{x}_1 \\[-0.15cm] \vdots \\ \dot{x}_N \end{bmatrix} =
    \begin{bmatrix} f_1(x_1,\dots,x_N) \\[-0.15cm] \vdots \\ f_N(x_1,\dots,x_N) \end{bmatrix} \hspace{-0.1cm},
\end{equation*}
where the time-scale separation in \eqref{eq:sysN} is recovered for $\tau_i=\frac{t}{\epsilon_i}$ in the \textit{singular limit} $\frac{\epsilon_{i+1}}{\epsilon_i} \to 0, \; \forall i=1,\dots,N-1$.

\subsection{Steady states, sensitivities and total derivatives}

\rr{Assume that every subsystem $\Sigma_i$ in \eqref{eq:sysN} has isolated steady states, and that the implicit function theorem \cite{krantz2012implicit} can be to guarantee the existence of steady state maps $x_i^s(x_1,\dots,x_{i-1})$ recursively from fast to slow subsystems $\Sigma_i$:} First, for some $x_N$ such that $0 = f_N(x_1,\dots,x_N)$ for the fastest $\Sigma_N$, the implicit function theorem guarantees the local existence of $x_N^s(x_1,\dots,x_{N-1})$. Under time-scale separation, the next system $\Sigma_{N-1}$ has reduced-order dynamics $\dot{x}_{N-1} = f^r_{N-1}(x_1,\dots,x_{N-1}) := f_{N-1}(x_1,\dots,x_{N-1},x_N^s(x_1,\dots,x_{N-1}))$, which allow now to define the steady state map $x_{N-1}^s(x_1,\dots,x_{N-2})$, and recursively $x_{i}^s(x_1,\dots,x_{i-1})$, which depend only on the states of the slower systems $x_1,\dots,x_{i-1}$. To ease the notation, from now on we will use $x_i^s$ instead of $x_{i}^s(x_1,\dots,x_{i-1})$ to denote the steady state map and its dependencies. Similarly, we use $(x_1,\dots, x_i,x_{i+1}^s,\dots,x_N^s)$ to denote that $x_{i+1}^s,\dots,x_N^s$ are all at steady-state for given values $x_1,\dots,x_i$. 
As in \eqref{eq:sys2ideal}, the steady state maps allow to define the reduced-order dynamics:

\begin{equation}\label{eq:sysNideal}
    \frac{dx_i}{d\tau_i} = f^r_i(x_1,\dots,x_i):= f_i(x_1,\dots, x_{i},x_{i+1}^s,\dots,x_N^s)
\end{equation}

Under sufficient regularity, the implicit function theorem gives the sensitivity of each steady state $x_i^s$ with respect to any $x_j$ for $j<i$:
\begin{equation}\label{eq:sysNsensss}\begin{array}{l}
     \nabla_{x_j}x_i^s(x_1,\dots,x_{i-1})  \\
     = - (\nabla_{x_i}f^r_i(x_1,\dots,x_{i-1},x_i^s))^{-1} \nabla_{x_j}f^r_i(x_1,\dots,x_{i-1},x_i^s) 
\end{array}
\end{equation}

Now the concepts of extended sensitivities \eqref{eq:sys2sens} and extended total derivatives \eqref{eq:bilsenstotderext} can be used to define compact analytical expressions for these sensitivities $\nabla_{x_j}x_i^{s}$ and an extension for a general point $(x_1, \allowbreak \dots, \allowbreak x_N)$. We define the extended total derivatives and sensitivities for every $i$ and $j$ recursively from $N$ to $1$: 
\begin{equation}\label{eq:sysNtotder}\begin{array}{rl}
D_{x_N}f_i := & \nabla_{x_N}f_i \\
D_{x_j}f_i := & \nabla_{x_j}f_i + \sum_{k=\max(i,j)+1}^N D_{x_k}f_i S_{x_{j}}^{x_{k}} \\
S_{x_j}^{x_i} := & - (D_{x_i}f_i)^{-1}D_{x_j}f_i, \; \forall j<i,
\end{array}
\end{equation}
where for clarity we omit the evaluation at $(x_1, \allowbreak \dots, \allowbreak x_N)$. These extended total derivatives of $f_i$ with respect to $x_j$ take into account the dependency of each intermediate $x_{k}^s$, $k \geq \max(i,j)+1$, with respect to $x_j$. When restricted, they coincide with the total derivatives \eqref{eq:bilsenstotderext} of $f_i^r(x_1,\dots,x_i)$, and the steady-state sensitivities \eqref{eq:sysNsensss}:
\begin{equation}\label{eq:sysNidealjac}\begin{array}{rl}
    D_{x_j}f_i(x_1, \dots, x_N)_{|_{(x_1,\dots,x_i,x_{i+1}^s,\dots,x_N^s)}} \hspace{-0.05cm} = & \hspace{-0.05cm} \nabla_{x_j}f^r_i(x_1,\dots,x_i) \\
    S_{x_j}^{x_i}(x_1, \dots, x_N)_{|_{(x_1,\dots,x_i,x_{i+1}^s,\dots,x_N^s)}} \hspace{-0.05cm} = & \hspace{-0.05cm} \nabla_{x_j}x_i^s(x_1,\dots,x_{i-1})
    \end{array}
\end{equation}

To guarantee that the implicit function theorem is applicable for every $f_i(x_1,\dots, x_{i},x_{i+1}^s,\dots,x_N^s)=0$, and thus ensure that steady states $x_i^s$, reduced-order dynamics \eqref{eq:sysNideal}, sensitivities \eqref{eq:sys2sensss}, and extended sensitivities and total derivatives \eqref{eq:sysNtotder} are well-defined around steady states, \rr{we formalize the assumptions made in this section in the following one:
\begin{ass}[Extension of Assumption~\ref{ass:sys2sensinv}] \label{ass:sysNsensinv}
For all $i>1$, the vector fields $f_i(\cdot),f_i^r(\cdot)$ are continuously differentiable, and the reduced-order systems \eqref{eq:sysNideal} have isolated steady-states $x_i^s$, where the partial derivatives $\nabla_{x_i}f_i^r(x_1,\dots,x_{i-1},x_i^s)$ are invertible. 
\end{ass}
This assumption also implies that the total derivatives $D_{x_i}f_i(\cdot)$ in \eqref{eq:sysNtotder} are invertible for $(x_i,\cdots,x_N)$ in a neighborhood of $(x_i^s,\cdots,x_N^s)$.}
Note that in contrast to Assumption~\ref{ass:sys2sensinv}, Assumption~\ref{ass:sysNsensinv} relaxes the need for a single steady-state in each subsystem, since there may exist multiple steady state $x_i^s$ for every subsystem $\Sigma_i$ given the values values $x_1,\dots,x_{i-1}$, i.e., the set of steady states
$$\mathcal{X}_i^s(x_1,\dots,x_{i-1}) := \{x_i | f_i(x_1,\dots,x_i,x_{i+1}^s,\dots,x_N^s)=0\}$$ 
is not necessarily a singleton. However, given the locally invertible $D_{x_i}f_i(\cdot)$, steady states will be isolated points. Moreover, given the implicit function theorem, the extended sensitivity $S_{x_j}^{x_i}(x_1,\dots,x_{i-1},x_i^s,\dots,x_N^s)$ \eqref{eq:sysNtotder} evaluated at each $x_i^s \in \mathcal{X}_i^s$ gives the actual sensitivity of each $x_i^s$ with respect to $x_1,\dots,x_{i-1}$. Therefore, Theorem~\ref{thm:sysNpredsensstab} presented later will allow to test the local stability of every combination of steady states $(x_1^s,\dots,x_N^s)$, where $x_i^s \in \mathcal{X}_i^s$ for every $i$.

\subsection{Sensitivity-conditioning for multiple time scales}

With the previously defined extended total derivatives and sensitivities \eqref{eq:sysNtotder}, the sensitivity-conditioning system \eqref{eq:sys2predsens} from Section~\ref{sec:predsens2sys} can be extended \rr{and applied to the multiple-time-scales system \eqref{eq:sysN}}:

\begin{equation}\label{eq:sysNpredsens}\begin{array}{l}
    \begin{bmatrix} I & 0 & 0 & 0 \\ 
    -S_{x_1}^{x_2} & I & 0 & 0 \\[-0.15cm]
    \vdots & \ddots & I & 0 \\
    -S_{x_1}^{x_N} & \cdots & -S_{x_{N-1}}^{x_N} & I \end{bmatrix} \hspace{-0.2cm}
     \begin{bmatrix} \dot{x}_1  \\[-0.15cm] \vdots \\ \dot{x}_N \end{bmatrix} 
     \hspace{-0.1cm} = \hspace{-0.1cm} \begin{bmatrix} f_1(x_1,  \dots,  x_N) \\[-0.15cm] \vdots \\ f_N(x_1, \dots, x_N) \end{bmatrix} \hspace{-0.1cm} ,
\end{array}
\end{equation}
where for clarity we omit the evaluation at point $(x_1, \allowbreak \dots, \allowbreak x_N)$ in the conditioning matrix $M$. In the equivalent expression
\begin{equation}\label{eq:sysNpredsenssysi}
\begin{array}{rl}
     \Sigma_i: \;\dot{x}_i = & f_i(x_1,\dots,x_N) + \sum_{j=1}^{i-1} S_{x_j}^{x_i}(x_1,\dots,x_N) \dot{x}_j, \\
\end{array}
\end{equation}
the extra terms $\sum_{j=1}^{i-1} S_{x_j}^{x_i}(x_1,\dots,x_N)\dot{x}_j$ play again the role of predicting and anticipating the changes of steady states $x_i^s$ due to the slower dynamics $\dot{x}_j, \; j <i$. \rr{Note that our approach \eqref{eq:sysNpredsens} does not requires any specific cascaded structure.}

As in Assumption~\ref{ass:sys2loclip}, to guarantee local existence and uniqueness \cite[Thm.~3.1]{khalil2002nonlinear} of a solution for \eqref{eq:sysNpredsens}, we assume:

\begin{ass}\label{ass:sysNloclip}
The vector fields in \eqref{eq:sysNpredsens},\eqref{eq:sysNpredsenssysi} are locally Lipschitz continuous.
\end{ass}

Then, the statements in Proposition~\ref{prop:sys2stab1} to \ref{prop:sys2stab3} for the two systems sensitivity-conditioning \eqref{eq:sys2predsens}, can be extended to the multiple time-scale case sensitivity-conditioning \eqref{eq:sysNpredsens}:

\begin{thm}[Extension of Propositions~\ref{prop:sys2stab1},\ref{prop:sys2stab2} and Corollary~\ref{cor:sys2stab3}]\label{thm:sysNpredsensstab} 
\
\begin{enumerate}[leftmargin=*]
\item\label{it:remain} \rr{\textbf{Positive invariance}: 
For some $p \geq 1$, given the dynamics \eqref{eq:sysNpredsenssysi} of $\dot{x}_i$ for $i \geq p$, initialized at time $t_0$: $x_i(t_0)=x_i^s(x_1(t_0),\dots,x_{i-1}(t_0))$. Then, $x_i(t)=x_i^s(x_1(t),\dots,x_{i-1}(t)) \; \forall i \geq p$ is the unique solution on the open domain of existence.
}

\item\label{it:locstab} \textbf{Local stability}: At a steady state $(x_1^s,\dots,x_N^s)$ the Jacobian $J$ of \eqref{eq:sysNpredsens} satisfies:
\begin{equation}\label{eq:sysNpredsensjac}\begin{array}{rl}
    J \hspace{-0.1cm} \sim \hspace{-0.15cm} & \hspace{-0.1cm}
    \begin{bmatrix} \nabla_{x_1}f_1^r(x_1^s) & \star & \star & \star \\
    0 & \nabla_{x_2}f_2^r(x_1^s,x_2^s) & \star & \star \\[-0.15cm]
    0 & 0 & \ddots & \star \\
    0 & 0 & 0 & \nabla_{x_N}f_N(x_1^s,\dots,x_N^s)
    \end{bmatrix} \hspace{-0.15cm} ,
\end{array}
\end{equation}
where $\star$ are irrelevant terms.

\item\label{it:globstab} \textbf{Global exponential stability}:
Assume that the vector fields $f_i(\cdot)$ are globally Lipschitz continuous, that the total derivatives $D_{x_i}f_i(\cdot)$ in \eqref{eq:sysNtotder} are globally invertible, and that there exists positive definite matrices $P_i \succ 0$ and $\eta_i$ such that for all $(x_1,\dots,x_N)$ the following contraction condition holds: 
\footnote{This would correspond to the condition ${P_i \nabla_{x_i}f_i^r + (\nabla_{x_i}f_i^r)^T P_i^T} \preceq - \eta_i P_i$ in {Corollary}~\ref{cor:sys2stab3} for the two time-scale system. However, here it is required to hold for an extended number of points $(x_1,\dots,x_N)$, not just for only $(x_1,\dots,x_{i},x_{i+1}^s,\dots,x_N^s)$, hence it is more strict.} 
\begin{equation}\label{eq:sysNcondstab}
{P_i D_{x_i}f_i+D_{x_i}f_i^T P_i^T} \preceq - \eta_i P_i.
\end{equation}
Then there exists a unique steady state $(x_1^s,\dots,x_N^s)$, which is a globally exponentially stable steady state of the sensitivity-conditioning approach \eqref{eq:sysNpredsens}.
\end{enumerate}
\end{thm}

\begin{proof}
See Appendix~\ref{app:thmproof} 
\end{proof}

\begin{cor}[of Theorem~\ref{thm:sysNpredsensstab}.\ref{it:locstab}, extending Corollary~\ref{cor:sys2corstab}]
The Jacobians \eqref{eq:sysNidealjac} and \eqref{eq:sysNpredsensjac}, for the systems under time-scale separation \eqref{eq:sysNideal} and the sensitivity-conditioning \eqref{eq:sysNpredsens}, respectively, have the same eigenvalues, and thus the same local stability. 
\end{cor}

The multiple time-scale sensitivity-conditioning \eqref{eq:sysNpredsens} can also be generalized as \eqref{eq:sys2predsensacc}, to improve the performance. Theorem~\ref{thm:sysNpredsensstab} can then be extended as in Proposition~\ref{prop:sys2stab4} and \ref{prop:sys2rob}.

\subsection{Discrete-time sensitivity-conditioning}\label{sec:predsensdisc}
Here we show how the multiple time-scale sensitivity-conditioning \eqref{eq:sysNpredsens} can be extended to discrete-time systems, while preserving the local stability result in Theorem~\ref{thm:sysNpredsensstab}. Consider the discrete-time systems: 
\begin{equation*}
 \Sigma_i: x_{i}^{k+1} = x_{i}^k + f_i(x_{1}^k,\dots,x_{N}^k),
\end{equation*}
where $x_{i}^k$ denote the value of $x_i$ at time $t^k$. 
Under time-scale separation as in \eqref{eq:sysNideal}, 
we can represent each discrete-time subsystem $\Sigma_i$ in its own time scale $k_i$ with reduced-order dynamics:
\begin{equation}\label{eq:sysNidealdisc}\begin{array}{rl}
    x_{i}^{k_i+1} = &  x_{i}^{k_i} + f^r_i(x_{1}^{k_i},\dots,x_{i-1}^{k_i},x_{i}^{k_i}) \\
     := & x_{i}^{k_i} + f_i(x_{1}^{k_i},\dots,x_{i-1}^{k_i},x_{i}^{k_i},x_{i+1}^s,\dots,x_{N}^s) 
\end{array}
\end{equation}
where $x_{j}^{k_i+1}=x_{j}^{k_i}$ for $j<i$, and $x_{j}^s$ is the steady state of $x_{j}$ given $x_1^{k_i},\dots,x_{i-1}^{k_i},x_{i}^{k_i},x_{i+1}^s,\dots,x_{j-1}^s$, for $j>i$. Similar to \eqref{eq:sysNidealjac} the Jacobians of \eqref{eq:sysNidealdisc} are
\begin{equation*}\begin{array}{l}
     J_i = I + \nabla_{x_i} f^r_i(x_{1}^{s},\dots,x_{i}^{s}) =  I +D_{x_i} f_i (x_{1}^s,\dots,x_{N}^s)
\end{array}
\end{equation*}
and thus systems \eqref{eq:sysNidealdisc} are locally asymptotically stable if the eigenvalues $\lambda_i$ of $I + \nabla_{x_i} f^r_i$ satisfy $\abs{\lambda_i}<1$.


Using the same extended sensitivities $S_{x_j}^{x_i}$ as in \eqref{eq:sysNtotder}, the discrete-time version of the sensitivity-conditioning \eqref{eq:sysNpredsens} can be expressed as:
\begin{equation}\label{eq:sysNpredsensdisc}\begin{array}{l}
    \begin{bmatrix} I & 0 & 0 & 0 \\ 
    -S_{x_1}^{x_2} & I & 0 & 0 \\[-0.15cm]
    \vdots & \ddots & I & 0 \\
    -S_{x_1}^{x_N} & \cdots & -S_{x_{N-1}}^{x_N} & I \end{bmatrix} \hspace{-0.2cm}
     \begin{bmatrix} x_1^{k+1}-x_{1}^k  \\[-0.15cm] \vdots \\ x_N^{k+1}-x_{N}^k \end{bmatrix} 
     \hspace{-0.1cm} = \hspace{-0.1cm} \begin{bmatrix} f_1 \\[-0.15cm] \vdots \\ f_N \end{bmatrix} \hspace{-0.1cm},
\end{array}
\end{equation}
omitting for clarity the evaluation at the point $(x_1^k,\dots,x_N^k)$.

\begin{prop}\label{prop:sysNpredsensstabdisc} At the steady state $(x_1^s,\dots,x_N^s)$ the Jacobians of \eqref{eq:sysNidealdisc} and \eqref{eq:sysNpredsensdisc} have the same eigenvalues. Thus, \eqref{eq:sysNidealdisc} is locally exponentially stable if and only if \eqref{eq:sysNpredsensdisc} is so. Moreover, if for some $i$ the Jacobian $I + \nabla_{x_i} f^r_i$ has any eigenvalue with norm larger than one, then both \eqref{eq:sysNidealdisc} and \eqref{eq:sysNpredsensdisc} are unstable. 
\end{prop}

\begin{proof} 
The proof follows similar steps as the local stability one in Theorem~\ref{thm:sysNpredsensstab}:
After performing the same similarity transformations, the Jacobian $J$ of \eqref{eq:sysNpredsensdisc} satisfies
\begin{equation*}
\begin{array}{rl}
    J \sim & I +
    \begin{bmatrix} \nabla_{x_1}f_1^r(x_1^s) & \star & \star \\[-0.15cm]
    0 & \ddots & \star \\
    0 & 0 & \nabla_{x_N}f_N(x_1^s,\dots,x_N^s)
    \end{bmatrix}     
\end{array}
\end{equation*}
\end{proof}

\section{Conclusion and Outlook}\label{sec:conc}

In this work, we have presented the sensitivity-conditioning: an alternative design tool for interconnected systems, that uses a predictive feed-forward term to preserve the stability of the system analysed at different time scales. This approach does not introduce a lower threshold on the actual time-scale separation between subsystems, in contrast to the usual singular perturbation approach. Moreover, we have shown examples of control design problems and optimization algorithms where our approach can be directly applied and improves the performance compared to a time-scale separation approach.

We believe that the applicability of our approach is not limited to these examples, but has the potential to be used in many other applications. For example, for any nested algorithms (e.g. in optimization or adaptive control) this sensitivity-conditioning could be used to design faster algorithms avoiding the need of time-scale separation between nested iterations. This is particularly promising for cases where iterations are computationally expensive, even when they are simple to evaluate, for example in distributed algorithms with communication bottlenecks.

Several directions for future research remain open: Since the sensitivities employed in the proposed conditioning are heavily model-based, we have established input-to-state stability robustness analysis against model errors. Nonetheless, a more sophisticated robust performance guarantees would be desirable. Moreover, through this work we have considered only continuously differentiable vector fields driving the dynamics. Thus, it remains to be seen how this method could be extended to nondifferentiable cases, arising often in optimization.


\appendices

\section{Proof of Lemma~\ref{lem:propext}}\label{app:propextproof}
For clarity we omit the evaluation at a given $x$. Consider the singular value decomposition $\nabla_{x}f=U S V^T$. Let $\sigma_{\min}=\min_i S_{i,i}$ denote the minimum singular value of $\nabla_{x}f$, and $v_{\min},u_{\min}$ the columns of $V,U$ corresponding to $\sigma_{\min}$. Since $P\succ 0$ and $\normsz{v_{\min}}_2=\normsz{u_{\min}}_2=1$, we have  
\begin{equation*}\begin{array}{rl}
    0 & < \eta \lambda_{\min}(P) \leq \eta v_{\min}^T P v_{\min} \\
    & {\leq} - v_{\min}^T (P \nabla_{x}f + \nabla_{x}f^T P^T) v_{\min}  \\
     & = - 2\sigma_{\min} v_{\min}^T P u_{\min}  \\
     & \leq 2\sigma_{\min} \normsz{v_{\min}}_2 \normsz{u_{\min}}_2 \lambda_{\max}(P) = 2\sigma_{\min} \lambda_{\max}(P), 
\end{array}
\end{equation*}
where the third inequality is due to $P\nabla_{x}f(x) \hspace{-0.05cm} \allowbreak + \hspace{-0.05cm} \nabla_{x}f(x)^T P^T \hspace{-0.05cm} \preceq \hspace{-0.05cm} -\eta P$ $\forall x$. Thus, $\nabla_{x}f^{-1} = V S^{-1} U^T$ is well-defined and  
$$ \normsz{\nabla_{x}f^{-1}}_2 = \max_i \tfrac{1}{S_{i,i}} = \tfrac{1}{\sigma_{\min}} \leq \tfrac{2\lambda_{\max}(P) }{\eta\lambda_{\min}(P)}.$$
Next, with $x(\delta) = x^s + \delta(x-x^s)$ for $\delta \in [0,1]$, then 
\begin{equation*}\begin{array}{l}
     (x-x^s)^T P \big( f(x) - \overbrace{f(x^s)}^{=0} \big) + \big( f(x) - \overbrace{f(x^s)}^{=0} \big)^T P (x-x^s)   \\
     = \int_0^1 (x-x^s)^T \big( \underbrace{P \nabla_x f(x(\delta)) + \nabla_xf(x(\delta)^T) P}_{\preceq -\eta P, \rr{\text{ since } x(\delta) \in \mathcal{B}_{r}(x^s)}}\big) (x-x^s) d\delta \\[-0.05cm]
     {\leq} -\eta\normsz{x-x^s}_P^2 \leq 0.
\end{array}
\end{equation*}
{Hence, $2\normsz{P^{\frac{1}{2}}f(x)}_2 = 2\normsz{f(x)}_P \normsz{x-x^s}_P \geq \eta \normsz{x-x^s}_P^2$,} and
$\normsz{f(x)}_P \geq \frac{\eta}{2} \normsz{x-x^s}_P$.

\section{Proof of Theorem~\ref{thm:sysNpredsensstab}}\label{app:thmproof}

\ref{it:remain}) {For any $i \geq p$, consider the dynamics of $\dot{x}_i$ in \eqref{eq:sysNpredsenssysi}, initialized at time $t_0$ as $x_i(t_0)=x_i^s(x_1(t_0),\dots,x_{i-1}(t_0))$. Local existence and uniqueness of a solution is guaranteed by Assumption~\ref{ass:sysNloclip}. Then, $x_i(t)=x_i^s(x_1(t),\dots,x_{i-1}(t))$ is the unique solution on the open domain of existence, since the derivative
\begin{equation*}\begin{array}{rl}
     & \dot{x}_i(t) \\
     = & f_i(x_1(t),\dots,x_N(t)) + \sum_{j=1}^{i-1} S_{x_j}^{x_i}(x_1(t),\dots,x_N(t)) \dot{x}_j(t) \\
    = & \hcancel{$f_i(x_1(t),\dots,x_{i-1}(t),x_{i}^s(t),\dots,x_N^s(t))$}{0}{0}{0}{0} \\
    & + \sum_{j=1}^{i-1} \underbrace{S_{x_j}^{x_i}(x_1(t),\dots,x_{i-1}(t),x_{i}^s(t),\dots,x_N^s(t))}_{\overset{\eqref{eq:sysNidealjac}}{=}\nabla_{x_j} x_i^s(x_1(t),\dots,x_{i-1}(t))} \dot{x}_j(t) \\[-0.1cm]
    = & \frac{dx_i^s(x_1(t),\dots,x_{i-1}(t))}{dt}
\end{array}
\end{equation*}
and the initial conditions coincide for $t \geq t_0$ on the domain of existence.}

\ref{it:locstab}) 
By recursively using the expression in \eqref{eq:sysNpredsenssysi}, the sensitivity-conditioning system \eqref{eq:sysNpredsens} can be rewritten as
\begin{equation*}\begin{array}{rl}
     \begin{bmatrix} \dot{x}_1 \\[-0.25cm] \vdots \\[-0.15cm] \dot{x}_N \end{bmatrix} \hspace{-0.15cm} 
     = & \hspace{-0.1cm} 
     \left[\begin{smallmatrix} I & 0 \\ 
     \begin{bmatrix} S_{x_1}^{x_N} \cdots S_{x_{N-1}}^{x_N} \end{bmatrix}  & I
     \end{smallmatrix}\right] \hspace{-0.2cm}
     \begin{bmatrix} \dot{x}_1 \\[-0.25cm] \vdots \\[-0.15cm] \dot{x}_{N-1} \\ f_N \end{bmatrix} \\[-0.1cm]
    = & \hspace{-0.1cm}
     \left[\begin{smallmatrix} I & 0 \\ 
     \begin{bmatrix} S_{x_1}^{x_N} \cdots S_{x_{N-1}}^{x_N} \end{bmatrix} & I
    \end{smallmatrix}\right]
     \hspace{-0.1cm} \cdots \hspace{-0.1cm}
     \left[\begin{smallmatrix} I & 0 & 0 \\ 
     \begin{matrix}S_{x_1}^{x_2}\end{matrix}  & I & 0 \\
     0 & 0 & I \end{smallmatrix}\right] 
     \hspace{-0.2cm}
     \begin{bmatrix} f_1 \\[-0.2cm] \vdots \\[-0.1cm] f_N \end{bmatrix} \hspace{-0.1cm},
\end{array}
\end{equation*}
where the matrices $\begin{bmatrix} S_{x_1}^{x_{i}} \cdots S_{x_{i-1}}^{x_{i}} \end{bmatrix}$ are at the row block $i$. Since at a steady state $(x_1^s,\dots,x_N^s)$ we have $f_i(x_1^s,\dots,x_N^s)=0$, then the Jacobian of \eqref{eq:sysNpredsens} at $(x_1^s,\dots,x_N^s)$ simplifies to

\begin{equation*}\begin{array}{l}
     J= \left[\begin{smallmatrix} I & 0 \\ 
     \begin{bmatrix} S_{x_1}^{x_N} \cdots S_{x_{N-1}}^{x_N} \end{bmatrix} & I
     \end{smallmatrix}\right] \cdots
    \left[\begin{smallmatrix} I & 0 & 0 \\ 
    \begin{matrix}S_{x_1}^{x_2}\end{matrix} & I & 0 \\
     0 & 0 & I \end{smallmatrix}\right]
      \nabla_{x} f, 
\end{array}
\end{equation*}
where the $i,j$-block of $\nabla_x f$ is $(\nabla_{x} f)_{i,j} = \nabla_{x_j} f_i$. 

Since  $\left[\begin{smallmatrix} I & 0 & 0 \\ 
     -\left[\begin{smallmatrix} S_{x_1}^{x_{i}} \cdots S_{x_{i-1}}^{x_{i}} \end{smallmatrix}\right]  & I & 0 \\
     0 & 0 & I \end{smallmatrix}\right]^{-1} = 
     \left[\begin{smallmatrix} I & 0 & 0 \\ 
     \left[\begin{smallmatrix} S_{x_1}^{x_{i}} \cdots S_{x_{i-1}}^{x_{i}} \end{smallmatrix}\right]  & I & 0 \\
     0 & 0 & I \end{smallmatrix}\right] \; \forall i$,
multiplying $J$ on both sides by these matrices we can iteratively construct matrices similar to $J$:
\begin{equation*}
\begin{array}{rl}
J \sim \tilde{J}_N = 
    & \left[\begin{smallmatrix} I & 0 \\ 
     -\begin{bmatrix} S_{x_1}^{x_N} \cdots S_{x_{N-1}}^{x_N} \end{bmatrix} & I
     \end{smallmatrix}\right] \hspace{-0.1cm}
     J \hspace{-0.1cm}
     \left[\begin{smallmatrix} I & 0 \\ 
     \begin{bmatrix} S_{x_1}^{x_N} \cdots S_{x_{N-1}}^{x_N} \end{bmatrix} & I
     \end{smallmatrix}\right] \\[0.3cm]
= &\left[\begin{smallmatrix} I & 0 & 0 \\ 
     \begin{bmatrix} S_{x_1}^{x_{N-1}} \cdots  S_{x_{N-2}}^{x_{N-1}} \end{bmatrix} & I & 0 \\
     0 & 0 & I \end{smallmatrix}\right] \cdots 
     \left[\begin{smallmatrix} I & 0 & 0 \\ 
    \begin{matrix}S_{x_1}^{x_2}\end{matrix} & I & 0 \\
     0 & 0 & I \end{smallmatrix}\right] \\[0.4cm]
    & \left[\begin{smallmatrix} \begin{bmatrix} (\nabla_{x_j}f_i+\nabla_{x_N}f_i S_{x_j}^{x_N})_{i,j}  \end{bmatrix} & \star \\
    0  & D_{x_N}f_N \end{smallmatrix}\right] \hspace{-0.1cm}, 
\end{array}
\end{equation*}
where $\begin{bmatrix} (\nabla_{x_j}f_i)_{i,j}  \end{bmatrix}$ indicates a matrix with block elements $(\nabla_{x_j}f_i)_{i,j}$. Assume that for $l+1$ we have
\begin{equation*}\begin{array}{l}
    \tilde{J}_{l+1} = \\
    \left[\begin{smallmatrix} I & 0 & 0 \\ 
     \begin{bmatrix} S_{x_{1}}^{x_{l}} \cdots S_{x_{N-2}}^{x_{l}} \end{bmatrix} & I & 0 \\
     0 & 0 & I \end{smallmatrix}\right] \cdots 
     \left[\begin{smallmatrix} I & 0 & 0 \\ 
    \begin{matrix}S_{x_1}^{x_2}\end{matrix} & I & 0 \\
     0 & 0 & I \end{smallmatrix}\right] \\[0.4cm]
     \left[\begin{smallmatrix} 
     \begin{bmatrix} (\nabla_{x_j}f_i + \sum_{k=l+1}^N D_{x_k}f_i S_{x_{j}}^{x_{k}} )_{i,j}  \end{bmatrix} & \star & \star & \star \\
    0 & D_{x_{l+1}}f_{l+1} & \star & \star \\[-0.25cm]
    0 & 0 & \ddots & \star \\
    0  & 0 & 0 & D_{x_N}f_N \end{smallmatrix}\right] \hspace{-0.1cm}
\end{array}
\end{equation*}
Note that the terms in the last row and column blocks of $\begin{bmatrix} (\nabla_{x_j}f_i + \sum_{k=j+1}^N D_{x_k}f_i S_{x_{j}}^{x_{k}} )_{i,j} \end{bmatrix}$ are total derivatives: $D_{x_j}f_l = \nabla_{x_j}f_l + \sum_{k=l+1}^N D_{x_k}f_l S_{x_{j}}^{x_{k}} $ for $j \leq l$, and $D_{x_l}f_i = \nabla_{x_l}f_i + \sum_{k=l+1}^N D_{x_k}f_i S_{x_{l}}^{x_{k}}$ for $i \leq l$. Therefore,  
\begin{equation*}\begin{array}{rl}
\tilde{J}_{l} \hspace{-0.1cm} =
     & \hspace{-0.1cm}  \left[\begin{smallmatrix} I & 0 & 0 \\ 
     -\begin{bmatrix} S_{x_{1}}^{x_{l}} \cdots S_{x_{N-2}}^{x_{l}} \end{bmatrix} & I & 0 \\
     0 & 0 & I\end{smallmatrix}\right] \hspace{-0.15cm}
     \tilde{J}_{l+1} \hspace{-0.15cm}
      \left[\begin{smallmatrix} I & 0 & 0 \\ 
     \begin{bmatrix} S_{x_{1}}^{x_{l}} \cdots S_{x_{N-2}}^{x_{l}} \end{bmatrix} & I & 0 \\
     0 & 0 & I \end{smallmatrix}\right] \\[0.4cm]
= & \hspace{-0.1cm}  \left[\begin{smallmatrix} I & 0 & 0 \\ 
     \begin{bmatrix} S_{x_{1}}^{x_{l-1}} \cdots S_{x_{N-2}}^{x_{l-1}} \end{bmatrix} & I & 0 \\
     0 & 0 & I \end{smallmatrix}\right] \cdots  
     \left[\begin{smallmatrix} I & 0 & 0 \\ 
    \begin{matrix}S_{x_1}^{x_2}\end{matrix} & I & 0 \\
     0 & 0 & I \end{smallmatrix}\right] \\[0.4cm]
     & \hspace{-0.1cm}  
     \left[\begin{smallmatrix} \begin{bmatrix} (\nabla_{x_j}f_i + \sum_{k=l}^N D_{x_k}f_i S_{x_{j}}^{x_{k}} )_{i,j}  \end{bmatrix} & \star & \star & \star \\
    \begin{bmatrix} (D_{x_{j}}f_{l} + D_{x_{l}}f_{l}S_{x_{j}}^{x_{l}} )_j \end{bmatrix} &  D_{x_{l}}f_{l} & \star & \star \\[-0.4cm]
    0 & 0 & \ddots & \star \\
    0  & 0 & 0 & D_{x_N}f_N \end{smallmatrix}\right] \hspace{-0.1cm},
\end{array}
\end{equation*}
where $D_{x_{j}}f_{l} + D_{x_{l}}f_{l}S_{x_{j}}^{x_{l}}=0 \; \forall j<l$ by definition of $S_{x_{j}}^{x_{l}}$ in \eqref{eq:sysNtotder}. Finally, $D_{x_1}f_1=\nabla_{x_1}f_1 + \sum_{k=2}^N D_{x_k}f_1 S_{x_{1}}^{x_{k}}$ implies:
\begin{equation*}\label{eq:sysNpredsensjactildep}\begin{array}{rl}
J\sim \tilde{J}_N \sim \dots \sim \tilde{J}_2 = &
    \left[\begin{smallmatrix} D_{x_1}f_1 & \star & \star \\[-0.25cm]
    0 & \ddots & \star \\
    0 & 0 & D_{x_N}f_N
    \end{smallmatrix}\right] 
\end{array}
\end{equation*}
and at steady-state $D_{x_i}f_i=\nabla_{x_i}f_i^r \; \forall i$.

\ref{it:globstab}) There is a unique steady state as a consequence of the contraction condition \eqref{eq:sysNcondstab} and Lemma~\ref{lem:propext}. The extended total derivatives \eqref{eq:sysNtotder} can be rewritten as 
\begin{equation}\label{eq:totderS}
    \begin{bmatrix} D_{x_1}f_i & \cdots & D_{x_N}f_i \end{bmatrix} 
    S_i
     =\begin{bmatrix} \nabla_{x_1}f_i & \cdots & \nabla_{x_N}f_i \end{bmatrix},
\end{equation}
where $S_i$ is the matrix with sensitivities in \eqref{eq:sysNpredsens} truncated at $i$:
\begin{equation*}
    S_i=\left[\begin{array}{ccc|ccc@{}} 
     & I &  &  & 0 &  \\ \hline
    -S_{x_1}^{x_{i+1}} & \cdots & -S_{x_{i}}^{x_{i+1}} & I & 0 & 0 \\[-0.15cm]
    \vdots & & & \ddots & I & 0\\
    -S_{x_1}^{x_N} & \cdots & -S_{x_{i}}^{x_N} & \cdots & -S_{x_{N-1}}^{x_N} & I\\
    \end{array}\right]
\end{equation*}

As in the proof of Proposition~\ref{prop:sys2stab3}, let $L_{f_i}$ denote the Lipschitz constant of $f_i$. Given Lemma~\ref{lem:propext}, it can be certified that the extended sensitivities $S_{x_{i}}^{x_{j}}$ and total derivatives $D_{x_j}f_i$ are all bounded, so there exists $L_{D,f_i,x_j}>0$ such that $\normsz{D_{x_j}f_i}_2 \leq L_{D,f_i,x_j}$. Defining the Lyapunov functions $V_i(x_1,\dots,x_N)=\normsz{f_i(x_1,\dots,x_N)}_{P_i}^2$, then we have
    \begin{equation*}\begin{array}{rl}
        
        \dot{V}_i = & f_i^T P_i \begin{bmatrix} \nabla_{x_1}f_i & \cdots & \nabla_{x_N}f_i \end{bmatrix} 
        \begin{bmatrix} \dot{x}_1^T  & \cdots & \dot{x}_N^T \end{bmatrix}^T\\
        & + (\text{transpose terms}) \\
        
        \overset{\eqref{eq:totderS}}{=} 
        & f_i^T P_i \begin{bmatrix} D_{x_1}f_i & \cdots & D_{x_N}f_i \end{bmatrix} S_i
        \begin{bmatrix} \dot{x}_1^T  & \cdots & \dot{x}_N^T \end{bmatrix}^T\\
        & + (\text{transpose terms}) \\
        
        \overset{\eqref{eq:sysNpredsens}}{=} 
        & f_i^T P_i \begin{bmatrix} D_{x_1}f_i & \cdots & D_{x_N}f_i \end{bmatrix}
        \begin{bmatrix} \dot{x}_1^T  & \cdots & \dot{x}_i^T & f_{i+1}^T & \cdots & f_{N}^T \end{bmatrix}^T\\
        & + (\text{transpose terms}) \\

        \overset{\eqref{eq:sysNpredsenssysi}}{=} 
        & f_i^T P_i \big( \cancel{\sum_{j=1}^{i-1} D_{x_j}f_i \dot{x}_j} + D_{x_i}f_i \underbrace{(f_i+ \cancel{\sum_{j=1}^{i-1} S_{x_j}^{x_i} \dot{x}_j} )}_{\dot{x}_i} \\[-0.7cm]
        & + \sum_{j=i}^N D_{x_j}f_i \; f_j \big) \\
        & + (\text{transpose terms})  \\
        
        \overset{\eqref{eq:sysNtotder}}{=}  & \sum_{j=i}^N f_i^T P_i (D_{x_j}f_i) f_j+ (f_i^T P_i (D_{x_j}f_i) f_j)^T \\
        
        \overset{\eqref{eq:sysNcondstab}}{\leq} &  - \eta_i \normsz{f_i}_{P_i}^2 + 2 \sum_{j=i+1}^N \underbrace{L_{D,f_i,x_j} \tfrac{\normsz{P_i}_2}{\lambda_{\min}({P_j})}}_{\nu_{i,j}} \normsz{f_i}_{P_i}\normsz{f_j}_{P_j} 
        
    \end{array}
    \end{equation*}

Consider the parameters $\theta_i>0$, and the Lyapunov function $V(x_1,\dots,x_N) = \sum_i \theta_i V_i(x_1,\dots,x_N)$. Then we have
\begin{equation*}\label{eq:sysNpredsensM} \begin{array}{rl}
    \dot{V} \leq & -\begin{bmatrix} \normsz{f_N}_{P_N} & \cdots &  \normsz{f_1}_{P_1} \end{bmatrix} W
    \begin{bmatrix} \normsz{f_N}_{P_N} & \cdots &  \normsz{f_1}_{P_1} \end{bmatrix}^T \\[0.1cm]

    W \hspace{-0.1cm} = & \hspace{-0.1cm}
    \begin{bmatrix} \theta_N  \eta_N & \cdots & -\theta_i \nu_{i,N} & \cdots & -\theta_1 \nu_{1,N} \\[-0.15cm] 
    \vdots & \ddots & \vdots & & \vdots \\[-0.1cm]
    -\theta_i \nu_{i,N} & \cdots & \theta_i \eta_i & \cdots & -\theta_1 \nu_{1,i} \\[-0.15cm] 
    \vdots & & \vdots & \ddots & \vdots  \\[-0.1cm]
    -\theta_1 \nu_{1,N} & \cdots & -\theta_1 \nu_{1,i} & \cdots & \theta_1  \eta_1 \end{bmatrix}
\end{array}
\end{equation*}
Let $W_i$ denote the upper left minor $i$ of $W$. The first minor satisfies $W_1= \zeta_1 := \theta_N \eta_N > 0$. We continue by recursion: assume there exists $\zeta_{i-1}>0$ so that $W_{i-1} \succeq \zeta_{i-1} I$, then $W_i \succ 0$ if 
    \begin{equation*}\begin{array}{l}
        \theta_{i} \eta_i - \theta_i^2
        \begin{bmatrix} \nu_{i,N} \\[-0.15cm] \vdots \\[-0.15cm] \nu_{i,i+1} \end{bmatrix}^T 
       \hspace{-0.1cm} W_{i-1}^{-1} 
        \begin{bmatrix} \nu_{i,N} \\[-0.15cm] \vdots \\[-0.15cm] \nu_{i,i+1} \end{bmatrix} \\[0.5cm]
        
        \geq \theta_i \eta_i - \frac{\theta_i^2}{\zeta_{i-1}}\sum_{j=i+1}^N \nu_{i,j}^2 > 0,
    \end{array}
    \end{equation*}
    or equivalently $\theta_i < \tfrac{\eta_i\zeta_{i-1}}{\sum_{j=i+1}^N \nu_{i,j}^2},$
which provides a recursive set of conditions on $\theta_i$ so that all minors $W_i$ are positive definite. In the end, $W_N=W \succ 0$, so there exists $\zeta_{N}>0$ such that $W \succeq \zeta_{N} I$ and
\begin{equation*}\begin{array}{c}
    \dot{V} \leq -\zeta_{N} \sum_{i=1}^N V_i \leq -\frac{\zeta_{N}}{\max_i \theta_i} V.
\end{array}
\end{equation*}




\bibliographystyle{IEEEtran}
\bibliography{ifacconf}

\end{document}